\documentclass [11pt,oneside]{amsart}

\reversemarginpar 
\usepackage{caption}

\usepackage{amssymb} \usepackage{amsfonts} \usepackage{amsmath}
\usepackage{amsthm} \usepackage{epsfig}

\usepackage{amssymb}
\usepackage{amsfonts}
\usepackage{amsmath}
\usepackage{amsthm}
\usepackage[english]{babel}
\usepackage{amscd,amsthm}
\usepackage{graphics}

\newcommand\noopsort[1]{}

\newtheorem{lemma}{Lemma}[section]
\newtheorem{teo}[lemma]{Theorem}
\newtheorem{prop}[lemma]{Proposition}
\newtheorem{cor}[lemma]{Corollary}

\theoremstyle{definition}
\newtheorem{defn}[lemma]{Definition}

\newtheorem{rem}[lemma]{Remark}

\newcommand{\matN}{\ensuremath {\mathbb{N}}}
\newcommand{\matR}{\ensuremath {\mathbb{R}}}

\newcommand{\matZ}{\ensuremath {\mathbb{Z}}}
\newcommand{\matC}{\ensuremath {\mathbb{C}}}

\newcommand{\matH}{\ensuremath {\mathbb{H}}}

\newcommand{\calI}{\ensuremath {\mathcal{I}}}

\newcommand{\calS}{\ensuremath {\mathcal{S}}}
\newcommand{\calM}{\ensuremath {\mathcal{M}}}

\newcommand{\calD}{\ensuremath {\mathcal{D}}}

\newcommand{\calT}{\ensuremath {\mathcal{T}}}

\newcommand{\tr}{{\rm tr}}

\newcommand{\finedimo}{{\hfill\hbox{$\square$}\vspace{2pt}}}
\newcommand{\dimostraz}{\vspace{2pt}\noindent\textsc{Proof:}\ }

\newcommand{\diag}{\SelectTips{cm}{} \xymatrix@1}

\def\diag{\SelectTips{cm}{} \xymatrix@1}
\newcommand{\tetra}{\Delta}

\newcommand{\TV}{{\rm TV}}

\newcommand{\nota} [1] {\caption{\footnotesize{#1}}}

\author{Roberto Frigerio}

\thanks{Supported by the INTAS project
``CalcoMet-GT'' 03-51-3663}

\title{{Similar fillings and isolation of cusps} 
{of hyperbolic $3$-manifolds}}

\keywords{Dehn filling, geodesic boundary, truncated tetrahedron, 
Kojima decomposition,
commensurability}

\address{Dipartimento di Matematica \\%
Universit\`a di Pisa \\%
Largo B.~Pontecorvo 5 \\%
56127 Pisa, Italy%
}

\email{frigerio@mail.dm.unipi.it}

\subjclass[2000]{57M50 (primary), 58H15 (secondary).}

\begin{document}

\begin{abstract}
In this paper we deepen the analysis
of certain classes $\calM_{g,k}$
of hyperbolic $3$-manifolds
that were introduced in a previous work by B.~Martelli,
C.~Petronio and the author.
Each element of $\calM_{g,k}$ is an oriented
complete finite-volume
hyperbolic $3$-manifold with compact connected geodesic boundary
of genus $g$ and $k$ cusps.
We study small deformations of the complete 
hyperbolic structure of manifolds in $\calM_{g,k}$
via a close analysis of their geodesic triangulations.
We
prove that several elements in $\calM_{g,k}$
admit non-homeomorphic hyperbolic Dehn fillings sharing the same
volume, homology, cusp volume, cusp shape,
Heegaard genus, complex length of the shortest geodesic, length
of the shortest return path, and Turaev-Viro invariants.
Manifolds which share all these invariants are called \emph{geometrically
similar}, and were first studied by C.~D.~Hodgson, R.~G.~Meyerhoff
and J.~R.~Weeks. The examples of geometrically similar manifolds
they described are commensurable with each
other.
We show here that many elements in  $\calM_{g,k}$ 
admit non-commensurable geometrically similar Dehn fillings.

The notion of 
\emph{geometric isolation} for cusps in a hyperbolic $3$-manifold
was introduced by W.~D.~Neumann and A.~W.~Reid 
and studied by D.~Calegary, who 
provided explanations 
for all the previously known examples of isolation
phenomena. We show here that
the cusps of any manifold $M\in \calM_{g,k}$ are 
geometrically isolated from each other.
Apparently, isolation of cusps in our examples arises
for different reasons from those described by Calegari.

We also show that
any element in $\calM_{g,k}$ admits an infinite
family of hyperbolic Dehn fillings inducing non-trivial deformations
of the hyperbolic structure on the geodesic boundary. 

\end{abstract}

\maketitle

\noindent
Let $N$ be an oriented complete finite-volume hyperbolic $3$-manifold with
compact geodesic boundary.
Mostow-Prasad's rigidity Theorem
implies that the space (of homotopy classes)
of complete finite-volume structures
supported by $N$ reduces to a single point, so
non-trivial deformations of the complete
structure can give rise only to incomplete metrics.
It is a well-known fact that 
such deformations are closely related 
to the geometry of manifolds which can be obtained 
from $N$ via \emph{Dehn filling}, as we are now going to explain.

A \emph{slope} 
on a torus is an isotopy class of simple unoriented closed curves.
Let $X$ be an oriented $3$-manifold with boundary tori 
$T_1,\ldots,T_k$ and let $V_1,\ldots,V_h$ be solid tori, $h\leqslant k$. 
Let $s_i$ be a slope on $T_i$ for $i=1,\ldots,h$ and choose
an attaching homeomorphism $\varphi_i:\partial V_i\to T_i$
taking a meridian of $V_i$ onto a loop representing $s_i$.
Set $\Phi=(\varphi_1,\ldots,\varphi_h)$ and  
$X (s_1,\ldots,s_h)=X\bigcup_\Phi (V_1\cup\ldots\cup V_h)$.
We say that $X(s_1,\ldots,s_h)$ is obtained by \emph{Dehn filling}
$X$ along the $s_i$'s. It is easily seen that
$X (s_1,\ldots,s_h)$ is a $3$-manifold whose homeomorphism
type depends solely on the $s_i$'s.
Also observe that the orientation of $X$ naturally induces an orientation
also on $X(s_1,\ldots,s_h)$.

A complete finite-volume hyperbolic $N$
admits a natural compactification obtained by adding
some boundary tori. 
Thus, up to identifying $N$ with its compactification, it does make sense 
to consider 
the Dehn fillings of $N$.
A crucial fact is that the metric completions
of many small deformations of the complete metric of $N$ actually
define complete hyperbolic structures on manifolds obtained by 
Dehn filling $N$. This phenomenon is at the heart of the proof of
Thurston's hyperbolic Dehn filling Theorem~\cite{Thu:bibbia}, which
states that ``almost all'' the Dehn fillings
of a cusped $3$-manifold support a complete finite-volume 
hyperbolic metric.

This paper is devoted to the description of certain classes $\calM_{g,k}$
of cusped hyperbolic $3$-manifolds with geodesic boundary.
Such classes were first introduced in~\cite{FriMaPe3}.
We concentrate here on 
describing how 
partially truncated triangulations can be employed in order
to study the 
Dehn fillings of manifolds in $\calM_{g,k}$.

\section{Preliminaries and statements}
All the manifolds considered in this paper will be connected
and oriented. Let $\Delta$ denote the standard tetrahedron, and let
$\dot{\Delta}$ be $\Delta$ with its vertices removed.
An \emph{ideal triangulation} of a compact $3$-manifold $M$ with boundary
is a realization of the interior of $M$ as a gluing of a finite number of
copies of $\dot{\Delta}$, induced by a simplicial face-pairing of 
the corresponding $\Delta$'s.  
Let $\Sigma_g$ be the closed orientable surface of genus $g$. The following
result~\cite{FriMaPe3} 
motivates the definition of $\calM_{g,k}$.
\begin{prop} \label{first:prop}
An ideal triangulation of a manifold whose boundary is the union of $\Sigma_g$
and $k$ tori
contains  at least $g+k$ tetrahedra.
\end{prop}
\noindent
For all $g>k\geqslant 1$
we then define $\calM_{g,k}$ as follows:
\begin{eqnarray*}
\calM_{g,k} & = & 
\big\{ {\rm compact\ oriented\ manifolds}\ M\ {\rm having\ an\ ideal\ triangulation}\\ 
& & \quad {\rm with\ } g+k\ {\rm tetrahedra,\ and\ }
\partial M = \Sigma_g\sqcup\big(\mathop{\sqcup}\limits_{i=1}^k T_i\big)\ 
{\rm with}\ T_i\cong\Sigma_1\big\}.
\end{eqnarray*}

Let $N$ be a compact manifold with boundary.
When this does not create ambiguities, we will denote
by $N$ also the manifold obtained by removing
the boundary tori from the original $N$. 
Thus the natural compactification
of a hyperbolic manifold will be usually denoted by the same symbol 
denoting the manifold itself. We say that a numerical sequence $\big(a_n\big)_{n=1}^\infty$ has 
\emph{growth type $n^n$} if there
exist constants $C>c>0$ such that $n^{c\cdot n}
< a_n < n^{C\cdot n}$ for $n\gg 0$. 
The following results are taken from~\cite{FriMaPe3}.

\begin{teo}
Any element in $\calM_{g,k}$ admits a complete finite-volume
hyperbolic structure with geodesic boundary.
\end{teo}
  
\begin{teo} \label{growth:teo}
For all $g>k\geqslant 1$ we have
$\calM_{g,k}\neq\emptyset$. 
Moreover, 
for any fixed $k$ the sequence $\big(\#\calM_{g,k}\big)_{g=2}^\infty$ 
has growth type $g^g$.
\end{teo}

\subsection{Isolation of cusps}
Recall that if $N$ is a complete finite-volume hyperbolic $3$-manifold,
then every boundary torus of $N$ is naturally endowed
with a Euclidean structure,
defined up to similarity.
Neumann and Reid introduced in~\cite{NeuRei} the notion of 
\emph{geometric isolation} for cusps in a hyperbolic manifold:

\begin{defn}\label{geo:isol:defn}
Let $N$ be a complete finite-volume hyperbolic
$3$-manifold with (possibly empty) geodesic boundary and
cusps $C_1,\ldots,C_h,C_{h+1},\ldots,C_k$. We say that
$C_1,\ldots,C_h$ are
geometrically isolated from $C_{h+1},\ldots,C_k$
if any small deformation of the hyperbolic structure on $N$
induced by Dehn filling $C_{h+1},\ldots,C_k$ while keeping
$C_1,\ldots,C_h$ complete does not affect the Euclidean structure
at $C_1,\ldots,C_h$.
\end{defn}

Calegary described in~\cite{Cal:napoleon} different
strategies for constructing manifolds with isolated cusps, also
providing explanations 
for all the previously known examples of isolation
phenomena. In Section~\ref{isolation:sec} we show that
the cusps of any manifold $M\in \calM_{g,k}$ are 
geometrically isolated from each other:

\begin{teo}\label{isolation:teo}
Let $M\in\calM_{g,k}$ with cusps $C_1,\ldots,C_k$ and let
$h\leqslant k$. Then $C_1,\ldots,C_h$ are
geometrically isolated from $C_{h+1},\ldots,C_k$.
\end{teo}

Apparently, isolation of cusps in our examples arises
for different reasons from those described in~\cite{Cal:napoleon}.

\subsection{Non-isolation of the boundary}
The natural question if the geodesic boundary of an element
in $\calM_{g,k}$ is isolated form the cusps is also answered.
Examples of isolation of the geodesic boundary from cusps of hyperbolic
$3$-manifolds were provided
in~\cite{NeuRei,Fuj:kodai}. On the other hand, non-isolation
phenomena were described in~\cite{Fuj:kodai2,FujKoj:osaka}.
In Section~\ref{nonisolation:section}
we prove the following:
\begin{teo}\label{nonisol:intro:teo}
Let $M\in\calM_{g,k}$. Then
there exists an infinite set $\{N_i\}_{i\in\matN}$ of complete 
finite-volume hyperbolic $3$-manifolds with the following property:
each $N_i$ is obtained by Dehn filling $M$, and the 
hyperbolic surfaces $\partial M, \partial N_1,\ldots,\partial N_i,\ldots$
are pairwise non-isometric.
\end{teo}

\subsection{Some invariants of hyperbolic 
$3$-manifolds}
Let $N$ be a complete hyperbolic $3$-manifold
and take a closed geodesic $\ell\subset N$. 
Then a well-defined \emph{complex length}
$\matC L(\ell)\in\matC/2\pi i\matZ$ exists which can
be described as follows. The universal covering
$\widetilde{N}$ of $N$ is isometric to a convex polyhedron
in $\matH^3$ bounded by a countable number oh hyperbolic 
planes~\cite{Koj:proc}.
Choose an orientation on $\ell$, and realize $\widetilde{N}$ 
in $\matH^3\cong\matC\times (0,\infty)$ in such a way that
$\ell$ lifts in $\widetilde{N}$ 
to the oriented geodesic $\widetilde{\ell}$
with endpoints $0$ and $\infty$. Let $\gamma\in{\rm Aut}
(\widetilde{N})\subset {\rm Isom}^+ (\matH^3)$ 
be the element corresponding to the oriented
curve $\ell$ which leaves $\widetilde{\ell}$ invariant.
A complex number $a$ exists such that
\[
\gamma (z,t)= (a\cdot z, |a| \cdot t),\quad (z,t)
\in\matC\times (0,\infty).
\]
We set $\matC L(\ell)=\ln a\in\matC/2\pi i\matZ$. 
It is easily seen that this is a good definition,
\emph{i.~e.}~that $a$ only depends on the \emph{unoriented} curve
$\ell$, and that
the usual length of $\ell$ 
is equal to $\Re (\matC L(\ell))$.

If $N$ is complete finite-volume 
with compact geodesic boundary and $k$ cusps, the \emph{cusp shape}
of $N$ is the set of Euclidean structures (up to a scale factor)
induced on the boundary tori of 
${N}$. A \emph{regular horocusp neighbourhood}
for $N$ is a set $O_1\sqcup\ldots\sqcup O_k\subset N$, where $O_i$
is an open embedded horospherical neighbourhood of the $i$-th cusp of $N$, 
$O_i\cap O_j=\emptyset$ for $i\neq j$ and ${\rm vol}(O_i)={\rm vol}
(O_j)$ for $i,j\in\{1,\ldots,k\}$. The \emph{cusp volume} of $N$
is the volume of a maximal regular horocusp neighbourhood for $N$
(where this volume is intended to be $0$ if $N$ is compact).
A \emph{return path} in $N$ 
is a geodesic segment in $N$ intersecting $\partial N$ perpendicularly
in its endpoints. Since the boundary of $N$ is compact, it is easily
seen that there exists a (not necessarily unique) shortest
return path in $N$. 

If $N$ is a compact 3-manifold with $\partial N=\partial_0N\sqcup\partial_1N$,
one can define the \emph{Heegaard genus} of $(N,\partial_0N,\partial_1N)$ 
as the minimal
genus of a surface that splits $N$ as $C_0\sqcup C_1$, where $C_i$
is obtained by attaching $1$-handles on one side of a collar of $\partial_iN$.
Moreover, for any integer $r\geqslant 2$,
after fixing in $\matC$ a primitive $2r$-th root of unity,
a real-valued invariant $\TV_r(N)$
was defined by 
Turaev and Viro in~\cite{TurVir}.

\subsection{Similar fillings} 
We are now ready to give the following:
\begin{defn}\label{geo:similar:def}
Let $N, N'$ be complete finite-volume hyperbolic $3$-manifolds 
with 
geodesic boundary and the same number of cusps. 
We say that $N$ and $N'$ are \emph{geometrically similar}
if the following conditions hold:
\begin{itemize}
\item
$N$ and $N'$ share the same volume, the same cusp volume 
and the same cusp shape;
\item
The shortest return paths of $N$ and $N'$ have the same length;
\item
The shortest closed 
geodesics of $N$ and $N'$ have the same complex length;
\item
$H_1(N;\matZ)\cong H_1(N';\matZ)$;
\item
if 
$\Sigma$ (resp.~$\Sigma'$) is the geodesic boundary of $N$ (resp.~of
$N'$) and $T_1,\ldots,T_k$ (resp.~$T'_1,\ldots,T'_k$)
are the boundary tori of $N$ (resp.~of $N'$), then
the Heegaard genus of $(N,\Sigma,T_1\sqcup\ldots\sqcup T_k)$
is equal to the Heegaard genus of $(N',\Sigma',T'_1\sqcup\ldots\sqcup T'_k)$;
\item
$N$ and $N'$ have the same Turaev-Viro invariants;
\item
Manifolds obtained by sufficiently complicated Dehn fillings on $N$
can be paired to manifolds   
obtained by sufficiently complicated Dehn fillings on $N'$
in such a way that the elements in each pair 
share the same volume, first homology group, cusp volume, cusp shape, length
of the shortest return path, complex length of the shortest geodesic,
Heegaard genus 
and Turaev-Viro invariants.
\end{itemize}
\end{defn}

Geometrically 
similar hyperbolic $3$-manifolds were first studied in~\cite{HodMeyWee}, where
it was shown that the Whitehead link complement admits
an infinite sequence of pairs of non-homeomorphic
geometrically similar Dehn fillings
(the definition of geometric similarity introduced in~\cite{HodMeyWee}
is actually a bit different from ours, and regards cusped
manifolds without geodesic boundary).
The elements of any pair of geometrically similar manifolds described
in~\cite{HodMeyWee} are both obtained by filling one cusp 
of the Whitehead link complement, and they are commensurable with each
other.
We show here that if $M\in\calM_{g,k}$ is generic,
\emph{i.~e.}~if it does not admit too many isometries,
then we can construct
different geometrically similar manifolds by filling $M$
along slopes on any chosen set of cusps of $M$.
This allows us to prove the following:

\begin{teo}\label{nonhomeo:fil:teo}
For any $k>0$ there exist $g>k$ and an element $X_k\in\calM_{g,k}$
with boundary tori $T_1,\ldots,T_k$
having the following property. For each $i=1,\ldots,k$
there exists a finite set $\calS_i$ 
of slopes on $T_i$ such that if
$h\leqslant k$
and $s_i\notin\calS_i$ is
a slope on $T_i$,
then $X_k (s_1,\ldots,s_h)$ is hyperbolic and 
at least
$(k!\cdot 3^h)/(h!\cdot (k-h)!)$ pairwise non-homeomorphic hyperbolic
Dehn fillings of $X_k$ are 
geometrically similar to $X_k (s_1\,\ldots,s_h)$.
\end{teo}

Moreover, the geometrically similar manifolds we obtain are typically
non-commensurable with each other (however, examples are also provided 
of non-homeomorphic geometrically similar commensurable
Dehn fillings of a specific element of $\calM_{g,k}$).

\section{Triangulations and deformation space}\label{tria:section}

In order to construct a hyperbolic structure on a
manifold $M\in\calM_{g,k}$, we 
choose a suitable triangulation of $M$ and we solve the
corresponding hyperbolicity equations.
We recall that the \emph{valence} of an edge
in a triangulation is the number of tetrahedra 
incident to it (with multiplicity).
The following result is proved in~\cite{FriMaPe3}. 

\begin{prop}\label{forma:tria:prop}
Let $M\in\calM_{g,k}$ 
and suppose
that 
$\calT$ is an ideal triangulation of $M$ with $g+k$ tetrahedra. Then 
the following holds:
\begin{itemize}
\item For any $i=1,\ldots,k$ there are exactly two tetrahedra
of $\calT$ with $3$ vertices on $\Sigma_g$ and one on $T_i$;
the remaining $g-k$ tetrahedra have all $4$ vertices on $\Sigma_g$;
\item
$\calT$ has $k+1$ edges $e_0,\ldots,e_k$ such that $e_0$ 
has both its endpoints on $\Sigma_g$ and valence $6g$,
while $e_i$ connects
$\Sigma_g$ to $T_i$ and has valence $6$
for $i=1,\ldots,k$.
\end{itemize}
\end{prop}

\subsection{Geometric tetrahedra}
In order to construct
a hyperbolic structure on our manifold $M\in\calM_{g,k}$ we realize the tetrahedra of an ideal 
triangulation of $M$ as special geometric blocks in $\matH^3$ and then 
we require that the structures match under the gluings. To describe the blocks
to be used we need some definitions. 

A \emph{partially truncated tetrahedron} is a pair $(\tetra,\calI)$, where 
$\tetra$ is a tetrahedron and either  $\calI=\emptyset$ 
or $\calI=\{v\}$, where $v$ is a vertex of $\Delta$.
In the latter case we say that $v$ is the
\emph{ideal vertex} of $\Delta$. In the sequel we will always
refer to $\tetra$ itself as a partially truncated tetrahedron, tacitly implying
that $\calI$ is also fixed. The \emph{topological realization} $\tetra^{\!\ast}$ 
of  $\tetra$
is obtained by removing from $\tetra$ the ideal vertex, if $\calI\neq\emptyset$,
and small open
stars of the non-ideal vertices. 
We call \emph{lateral hexagon} and \emph{truncation 
triangle} 
the intersection of $\tetra^{\!\ast}$ respectively with a face of $\tetra$ 
and with the link in $\tetra$ of a non-ideal vertex. The edges of the 
truncation triangles, which also belong to the lateral hexagons, are called 
\emph{boundary edges}, and the other edges of $\tetra^{\!\ast}$ are called
\emph{internal edges}. If $\tetra$ has an ideal vertex, 
three lateral
hexagons of $\tetra^{\!\ast}$ are in fact pentagons with a vertex removed,
and they are called
\emph{exceptional lateral hexagons}.

A \emph{geometric realization} of $\tetra$ is an 
embedding of $\tetra^{\!\ast}$ in $\matH^3$ such that the truncation triangles
are geodesic triangles, the lateral hexagons are 
geodesic polygons with ideal vertices
corresponding to missing edges, and truncation 
triangles and lateral hexagons 
lie at right angles to each other.
The following theorem~\cite{Fuj:tokyo,FriPe} classifies isometry classes 
of geometric partially truncated tetrahedra.

\begin{teo}\label{moduli:teo}
Let $\Delta$ be a partially truncated tetrahedron and let 
$\Delta\!^{(1)}$ be the
set of edges of $\Delta$. The geometric realizations of $\Delta$ are
parameterized up to isometry 
by the dihedral angle assignements $\theta:\Delta\!^{(1)}\to(0,\pi)$ such that
for each vertex $v$ of $\Delta$, 
if $e_1,e_2,e_3$ are the edges that emanate from $v$, then
$\theta(e_1)+\theta(e_2)+\theta(e_3)$ is equal to $\pi$ for ideal $v$
and less than $\pi$ for non-ideal $v$.
\end{teo}

The following well-known hyperbolic trigonometry formulae will prove useful
later:
\begin{figure}
\begin{center}
\input{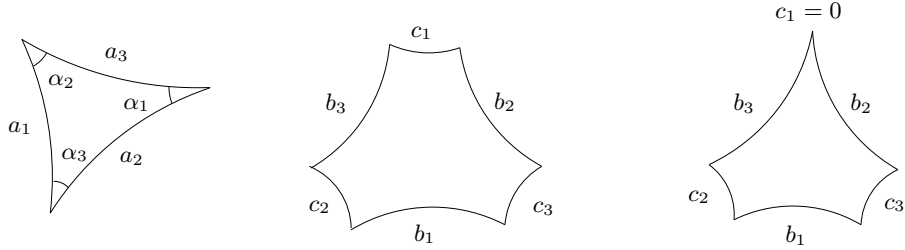}
\caption{{A triangle, a right-angled hexagon and a pentagon
with four right angles and an ideal vertex.}}\label{trigo:fig}
\end{center}
\end{figure}

\begin{lemma}\label{trigo:lemma}
With notation as in Fig.~\ref{trigo:fig} we have
\begin{eqnarray}
\cosh a_1={(\cos \alpha_2\cdot\cos \alpha_3+\cos\alpha_1)}/{(\sin
\alpha_2\cdot\sin\alpha_3)},\label{cos:rule}\\
\sinh a_1/\sin\alpha_1=\sinh a_2/\sin \alpha_2=\sinh a_3/\sin\alpha_3,
\label{sin:rule}\\
\cosh b_1={(\cosh c_2\cdot\cosh c_3+\cosh c_1)}/{(\sinh
c_2\cdot\sinh c_3)}\label{hexa:rule}.
\end{eqnarray}
\end{lemma}

\subsection{Hyperbolicity equations}\label{hyp:subsection}
Let $M$ be an element of $\calM_{g,k}$ and
$\calT$ be an ideal triangulation of $M$ with $g+k$
tetrahedra.
We try to give $M$ a hyperbolic structure
with geodesic boundary by looking for a geometric
realization $\theta$ of $\calT$
such that the structures of the tetrahedra
match under the gluings. 
In order to define a global hyperbolic structure on $M$,
the tetrahedra of $\calT$ must satisfy two obvious necessary 
conditions, which in fact are also sufficient. Namely, 
we should be able
to glue the lateral hexagons by isometries, and
we should have a total dihedral angle of $2\pi$ around each 
edge of the manifold.
The first condition ensures that the hyperbolic structure defined
by $\theta$ on the complement of the $2$-skeleton of $\calT$
extends to the complement of the $1$-skeleton.  Since any tetrahedron
of $\calT$ contains at most one ideal vertex,
the second one ensures that such a structure glues up without singularities 
also along the edges.

By Proposition~\ref{forma:tria:prop}, 
if we suppose $M$ to be hyperbolic and $\calT$
to be geometric (\emph{i.~e.}~to define a hyperbolic structure on
the whole of $M$), than the  
edges of the tetrahedra with all the vertices on $\Sigma_g$
should have all 
the same length. This would force the realizations of
the compact tetrahedra 
in $\calT$ to be regular and isometric to each other.
Moreover, all the finite internal edges of the 
tetrahedra with one vertex on
the boundary tori should also have the same length.

\begin{figure}
\begin{center}
\input{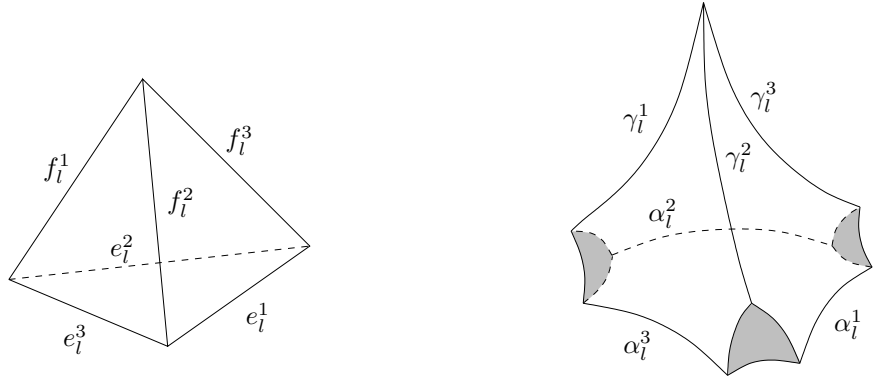}
\caption{{Dihedral angles along the edges of a non-compact 
tetrahedron of $\calT$.}}\label{par:fig}
\end{center}
\end{figure}

On each tetrahedron of $\calT$ we fix the orientation compatible
with the global orientation of $M$. As a result also the lateral hexagons
have a fixed orientation, which is reversed by the gluing maps.
Let us now fix some notation we will use extensively later on.
Let $T_1,\ldots, T_k$ be the boundary tori of $M$.
We denote by $\Delta_{2i-1},\Delta_{2i}$ the tetrahedra of $\calT$
incident to $T_i$ and by $F_{2i-1}^1,F_{2i-1}^2,F_{2i-1}^3,F_{2i}^1,
F_{2i}^2,F_{2i}^3$ the exceptional hexagons
of $\Delta_{2i-1},\Delta_{2i}$, in such a way that 
$F_{2i-1}^j$ is glued to $F_{2i}^j$
for $j=1,2,3$. For $l=1,\ldots,2k$
we also suppose that $F^1_l,F^2_l,F^3_l$ are positively arranged 
around the ideal
vertex of $\Delta_l$, and we call
$e_l^j$ the only finite internal edge of $F_l^j$, and $f_l^j$
the edge of $\Delta_l$ opposite to $e_l^j$. 
We now consider a geometric realization $\theta$
of the tetrahedra
of $\calT$ such that compact tetrahedra are regular and isometric 
to each other, and for $l=1,\ldots,2k$, $j=1,2,3$
we set $\alpha_l^j=\theta(e_l^j)$, and $\gamma_l^j=\theta
(f_l^j)$ (see Fig.~\ref{par:fig}). 
We set $\beta$ to be the
dihedral angle along the edges of the $g-k$ compact tetrahedra
of $\calT$. We denote by $L^\theta$ the length with respect to the
realization $\theta$.

\subsection{Consistency along the faces}
We first determine the conditions on dihedral angles
under which  
all the \emph{compact} lateral hexagons of the tetrahedra
in $\calT$ are regular and isometric to each other.
By equation~(\ref{hexa:rule}), this is equivalent to asking
that the lengths of all the boundary edges of all the compact  
lateral hexagons are equal to each other, and 
by equation~(\ref{cos:rule}),
this condition translates into the following set of equations:
\begin{equation}\label{length:edge:eq}
\frac{\cos\alpha_l^j\cdot\cos\alpha_{l}^{j+1}+\cos\gamma_{l}^{j+2}}
{\sin\alpha_l^j\cdot\sin\alpha_{l}^{j+1}}
=\frac{\cos^2 \beta+\cos\beta}{\sin^2\beta},\qquad l=1,\ldots,2k,\ j=1,2,3.
\end{equation}

\subsection{Exceptional hexagons}

\begin{figure}
\begin{center}
\input{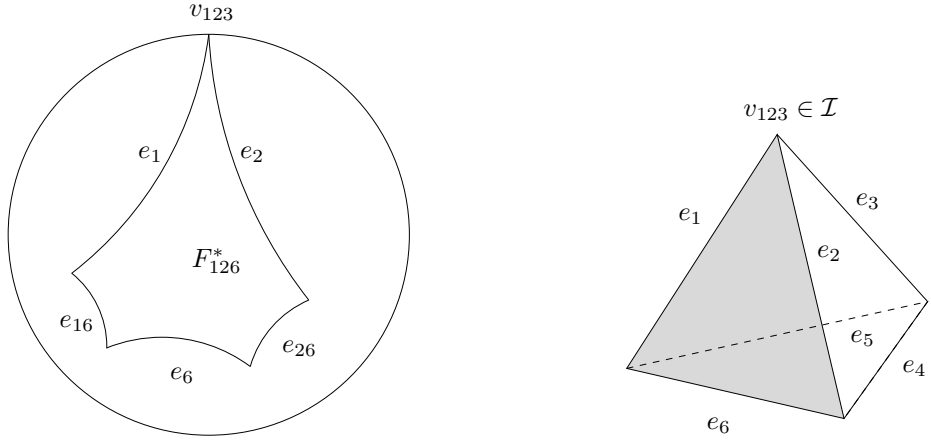}
\caption{{An exceptional hexagon.}}\label{special:hexagon:fig}
\end{center}\end{figure}

Let us consider an exceptional hexagon $F^*_{126}$ 
as in Fig.~\ref{special:hexagon:fig}, and  recall that the hexagon is
embedded in $\matH^3$ by $\theta$.
We consider the horospheres $O_1$ and
$O_2$ centred at $v_{123}$ and passing through the 
non-ideal ends of $e_1$ and $e_2$ respectively.
We define $\sigma^\theta(F_{126})$ to be
$\pm{\rm dist}(O_1,O_2)$, the sign being positive
if $e_2,v_{123},e_1$
are arranged positively on $\partial F^*_{126}$
and $O_1$ is contained in the horoball bounded
by $O_2$, or if $e_2,v_{123},e_1$ are arranged negatively
on $\partial F^*_{126}$ and $O_2$ is contained in the horoball
bounded by $O_1$, and 
negative otherwise. 
Let us denote by $e_{ij}$ the boundary edge joining $e_i$ with $e_j$.
From~\cite[Proposition 1.8]{Fri:preprint}
and equation~(\ref{sin:rule}) we deduce:
\begin{equation}\label{sigma:compute:eq}
\sigma^\theta(F_{126})=\ln \frac{\sinh L^\theta (e_{56})}
{\sinh L^\theta (e_{46})}+\ln \frac{\sin \theta (e_2)\cdot\sin\theta (e_5)}
{\sin\theta (e_1)\cdot\sin\theta (e_4)}.
\end{equation}

We now set $\ell^\theta(F_{126})=L^\theta (e_6)$.
The next proposition~\cite{Fri:preprint} shows that the functions 
$\sigma$ and $\ell$ provide a parameterization
of isometry classes of exceptional hexagons.

\begin{prop}\label{sigma:works:prop}
Let $F$ and $F'$ be paired exceptional lateral hexagons.
Their pairing can be realized by an isometry if and only
if $\sigma^\theta(F)+\sigma^\theta(F')=0$
and $\ell^\theta(F)=\ell^\theta(F')$.
\end{prop}

Recall now that we are considering 
the geometric realization of the triangulation $\calT$ of $M$
parameterized by the dihedral angles $\alpha_i^j, \gamma_i^j,\beta$,
$i=1,\ldots,2k,\, j=1,2,3$.
Under the assumption that
equations~(\ref{length:edge:eq}) are in force (\emph{i.~e.}~that
the boundary edges of the compact faces of $\calT$ have all the same length),
Proposition~\ref{sigma:works:prop} and equation~(\ref{sigma:compute:eq})
imply that the matching exceptional hexagons can be glued by isometries
if and only if for $i=1,\ldots,k$ we have:
\begin{equation}\label{sigma:eq}
\begin{array}{llll}
& \sin\alpha_{2i}^1\sin\alpha_{2i-1}^1\sin\gamma_{2i}^1\sin\gamma_{2i-1}^1&=&
\sin\alpha_{2i}^2\sin\alpha_{2i-1}^2\sin\gamma_{2i}^2\sin\gamma_{2i-1}^2\\
= & \sin\alpha_{2i}^3\sin\alpha_{2i-1}^3\sin\gamma_{2i}^3\sin\gamma_{2i-1}^3.
& &
\end{array}
\end{equation}

\subsection{Consistency around the edges}
Since $\gamma_l^1+\gamma_l^2+\gamma_l^3=\pi$
for $l=1,\ldots,2k$, the total angle along
any half-infinite edge of $\calT$ is automatically forced
to be equal to $2\pi$, so consistency around the edges 
translate into
the following equation only:
\begin{equation}\label{sum:eq}
6\cdot (g-k)\cdot \beta+\sum_{l=1}^{2k} (\alpha_l^1+\alpha_l^2+\alpha_l^3)
=2\pi.
\end{equation}
Any solution of 
\emph{consistency equations}~(\ref{length:edge:eq}),~(\ref{sigma:eq}),~(\ref{sum:eq}) 
defines a non-singular hyperbolic structure with geodesic boundary on $M$.

\subsection{Completeness equations}\label{compl:subsection}
For $i=1,\ldots,k$ let now $\mu_i,\lambda_i$ be the basis of
$H_1 (T_i;\matZ)$ which is defined as follows:
$\mu_i$ is the projection on $T_i$ of the edge
in the link of the ideal vertex of $\Delta_{2i-1}$
that joins $f_{2i-1}^1$ to $f_{2i-1}^2$;
$\lambda_i$ is the projection on $T_i$ of the edge
in the link of the ideal vertex of $\Delta_{2i}$
that joins $f_{2i}^3$ to $f_{2i}^2$.
A solution
\begin{equation}\label{xxxx:eq}
x =
 (\alpha_1^1,\alpha_1^2,\alpha_1^3,
\gamma_1^1,\gamma_1^2,\gamma_1^3,\ldots,
\alpha_{2k}^1,\alpha_{2k}^2,\alpha_{2k}^3,
\gamma_{2k}^1,\gamma_{2k}^2,\gamma_{2k}^3, \beta)\in\matR^{12k+1}
\end{equation}
of consistency equations naturally defines an ${\rm Aff} (\matC)$-structure
on $T_i$ (see \emph{e.~g.}~\cite{BenPet:book,Fri:tesi}).
We denote by 
$a_i (x)\in\matC$ (resp.~by $b_i (x)\in\matC$) the dilation component
of the holonomy of $\mu_i$ (resp.~of $\lambda_i$) corresponding
to the ${\rm Aff} (\matC)$-structure defined by $x$ on $T_i$.
It is well-known that the hyperbolic
structure defined by $x$ on $M$ induces a complete metric on
the $i$-th end of $M$ if and only if $a_i (x)=b_i (x)=1$.
Moreover, one can explicitly compute $a_i$ and $b_i$
in terms of the dihedral angles:
\begin{figure}
\begin{center}
\input{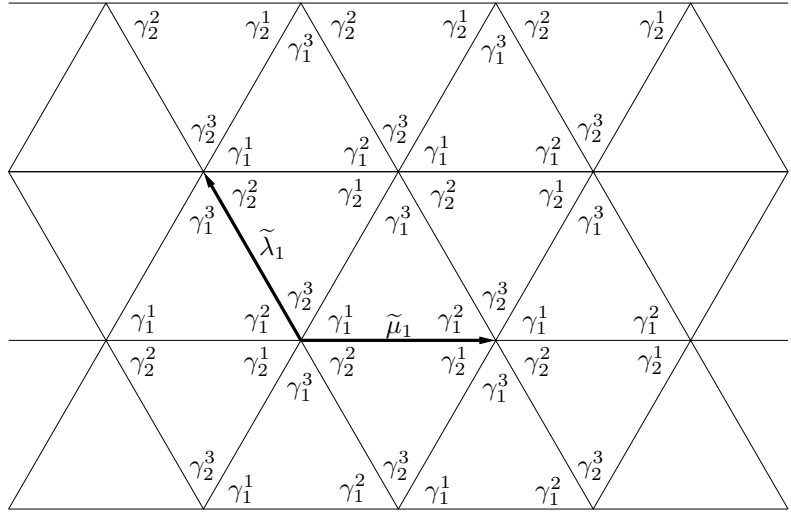}
\caption{{The triangulation of the universal covering of 
$T_1$.}}\label{tessel2:fig}
\end{center}
\end{figure}

\begin{teo}\label{computing:hol:teo}
We have
\[
\begin{array}{lll}
a_i (x)&=&((\sin\gamma_{2i-1}^1\sin\gamma_{2i}^2)/
(\sin\gamma_{2i-1}^2\sin\gamma_{2i}^1))\exp (i(\gamma_{2i-1}^3-\gamma_{2i}^3)),
\\
b_i (x)&=&((\sin\gamma_{2i-1}^2\sin\gamma_{2i}^3)/
(\sin\gamma_{2i-1}^3\sin\gamma_{2i}^2))\exp (i(\gamma_{2i-1}^1-\gamma_{2i}^1)).
\end{array}
\]
\end{teo}
\dimostraz
Apply~\cite[Proposition 1.14]{Fri:preprint} (see Fig.~\ref{tessel2:fig}).
\finedimo

\subsection{The complete solution}\label{compl2:subsection}
The following theorem~\cite{FriMaPe3} shows that a solution of
consistency and completeness equations always 
exists, and is as symmetric
as possible.

\begin{teo}\label{solution:teo}
There exist constants
$\overline{\alpha}_{g,k},\overline{\beta}_{g,k}
\in (0,\pi/3)$ such that the point
\[
x_0=
(\overline{\alpha}_{g,k},\overline{\alpha}_{g,k},
\overline{\alpha}_{g,k},
\pi/3,\pi/3,\pi/3,\ldots,
\overline{\alpha}_{g,k},\overline{\alpha}_{g,k},\overline{\alpha}_{g,k},
\pi/3,\pi/3,\pi/3, \overline{\beta}_{g,k})\in\matR^{12k+1}
\]
provides the \emph{unique}
solution of consistency and completeness equations for $\calT$.
\end{teo}
Thus 
the complete hyperbolic structure of $M$ induces on each boundary
torus the regular hexagonal Euclidean structure which is obtained
by gluing two Euclidean equilateral triangles.

\begin{lemma}\label{bgg:lemma}
We have
$\overline{\alpha}_{g,k}<\overline{\beta}_{g,k}<2\overline{\alpha}_{g,k}
\leqslant \pi/3$.
\end{lemma}
\dimostraz
From equation~(\ref{length:edge:eq})
we easily get
$\cos\overline{\beta}_{g,k}=(2+\cos 2\overline{\alpha}_{g,k})/3>\cos 2\overline{\alpha}_{g,k}$,
whence $\overline{\beta}_{g,k}<2\overline{\alpha}_{g,k}$. 
Moreover, since $\overline{\alpha}_{g,k}<\pi/3$
we have $2\cos^2\overline{\alpha}_{g,k}-3\cos\overline{\alpha}_{g,k}+1<0$, whence
$\cos\overline{\beta}_{g,k}=(2\cos^2\overline{\alpha}_{g,k}+1)/3<\cos\overline{\alpha}_{g,k}$
and $\overline{\alpha}_{g,k}<\overline{\beta}_{g,k}$. This inequality also
implies
$2\pi=6(g-k)\cdot\overline{\beta}_{g,k}+6k\cdot \overline{\alpha}_{g,k}\geqslant
6g \cdot \overline{\alpha}_{g,k}$, whence
$\overline{\alpha}_{g,k}\leqslant \pi/6$.
\finedimo

\begin{lemma}\label{embedcusp:lemma}
Let $\Delta^{\!\ast}\subset\matH^3$ be a
non-compact geometric tetrahedron
of the realization of $\calT$ 
parameterized by $x_0$ and let $v$ be the ideal vertex
of $\Delta^{\!\ast}$.
Then the horosphere 
centred at $v$ and tangent to the truncation triangles of
$\Delta^{\!\ast}$
does not intersect the lateral
hexagon opposite to $v$.
\end{lemma}
\dimostraz
For $\epsilon>0$ let $\Delta^{\!\ast}_\epsilon$ be the small deformation
of $\Delta^{\!\ast}$ having angle  $\pi/3-\epsilon$
along the edges emanating from $v_\epsilon$,
and $\overline{\alpha}_{g,k}$
along the other internal edges.
Let $l_\epsilon$ be the length of an internal edge emanating
from $v_\epsilon$ and $d_\epsilon$ be the distance between the truncation
triangle corresponding to $v_\epsilon$ and the opposite lateral hexagon.
By~\cite[Proposition~5.6]{FriPe}, Lemma~\ref{trigo:lemma} and 
Lemma~\ref{bgg:lemma}
we have 
\[
\lim_{\epsilon\to 0} (\cosh d_\epsilon/\cosh l_\epsilon)=\sqrt{4\cos^2
\overline{\alpha}_{g,k}-1}>1,
\]
whence the conclusion.
\finedimo

\subsection{Smoothness at the complete structure}\label{smooth:subsection}
From now on we denote
by $\Omega_{g,k}\subset\matR^{12k+1}$
the set of solutions of consistency equations for $\calT$
(it is clear that this set indeed depends only on $g$ and
$k$, and not on $\calT$).
If $x\in\Omega_{g,k}$ is as in equation~(\ref{xxxx:eq}),
we set
\begin{equation}\label{uv:eq}
\begin{array}{lllll}
u_i (x)&=& \ln a_i (x) & = & 
\ln ({(\sin\gamma_{2i-1}^1\sin\gamma_{2i}^2})/
{(\sin\gamma_{2i-1}^2\sin\gamma_{2i}^1)})+ i(\gamma_{2i-1}^3-\gamma_{2i}^3),
\\
v_i (x)&=& \ln b_i (x)&=&\ln ({(\sin\gamma_{2i-1}^2\sin\gamma_{2i}^3)}/
{(\sin\gamma_{2i-1}^3\sin\gamma_{2i}^2)})+ i(\gamma_{2i-1}^1-\gamma_{2i}^1).
\end{array}
\end{equation}
The following theorem is proved in~\cite{Fri:preprint} (see also~\cite{NeuZag}).

\begin{teo}\label{smoothness:teo}
Near $x_0$, the space $\Omega_{g,k}\subset\matR^{12k+1}$ is a smooth manifold
of real dimension $2k$, whose 
tangent space $T_{x_0}\Omega_{g,k}$ at $x_0$ is given 
by the solutions of the linearization of consistency
equations~(\ref{length:edge:eq}),~(\ref{sigma:eq}),~(\ref{sum:eq}).
Moreover,
there exists a small neighbourhood $U$ of $x_0$ in $\Omega_{g,k}$ 
with the following properties:
\begin{enumerate}
\item
For $x\in U$, we have $u_i (x)=0$ $\Leftrightarrow$ $v_i (x)=0$
$\Leftrightarrow$ the metric structure defined by $x$ is complete at
the $i$-th end of $M$;
\item
The map
$(u_1,\ldots,u_k):U\to\matC^k$ is a diffeomorphism
between $U$ and an open neighbourhood of $0$ in $\matC^k$.
\end{enumerate}
\end{teo}

Due to our choice of $\mu_i,\lambda_i$ we also have the following:

\begin{lemma}\label{neuzag:uv:lemma}
If $j\in\{1,\ldots,k\}$ and  
$\{y_n\}_{n\in\matN}\subset \Omega_{g,k}$ is a sequence with
$\lim_{n\to\infty} y_n=x_0$ and $u_j (y_n)\neq 0$ for every $n\in\matN$, then
$\lim_{n\to\infty} v_j (y_n)/u_j (y_n)=-1/2+i\sqrt{3}/2$.
\end{lemma}
\dimostraz See~\cite{NeuZag,Fri:tesi}.
\finedimo

\subsection{Dehn filling equations}\label{Dehn:subsection}
Let $U$ be a sufficiently small neighbourhood of $x_0$
in $\Omega_{g,k}$ and let $x\in U$. For $j=1,\ldots,k$, we define
the $j$-\emph{Dehn filling
coefficient} $(p_j (x),q_j (x))\in\matR^2\cup\{\infty\}$
as follows: if $u_j (x)=0$, then $(p_j (x),q_j (x))=\infty$;
otherwise, $p_j (x),q_j (x)$ are the unique real solutions
of the equation
\[
p_j (x) u_j (x)+q_j (x) v_j (x)=2\pi i.
\]
(Existence and uniqueness of such solutions
near $x_0$ can be easily deduced from Theorem~\ref{smoothness:teo}
and Lemma~\ref{neuzag:uv:lemma}.)
Let us set
\[
d=(d_1,\ldots,d_k):U\to \prod_{i=1}^k S^2,\quad
d_j (x)=(p_j(x),q_j(x))\in S^2=\matR^2\cup\{\infty\}.
\]
As a consequence of Theorem~\ref{smoothness:teo}
and Lemma~\ref{neuzag:uv:lemma}
we have the following:

\begin{teo}\label{essential:bis:teo}
If $U$ is small enough,
the map $d$ 
defines a
diffeomorphism onto
an open neighbourhood of
$(\infty,\ldots,\infty)$ in $S^2\times\cdots\times S^2$.
\end{teo}

For $x\in\Omega_{g,k}$ we denote by $M(x)$
the hyperbolic structure induced on $M$ by $x$, and by
$\widehat{M} (x)$ the metric completion of 
$M(x)$. We also set
\begin{eqnarray*}
I\Omega_{g,k}&= &\big\{x\in U\subset\Omega_{g,k}:
\ {\rm for}\ i=1,\ldots,k\ {\rm the}
\ i-{\rm th\ Dehn\ filling\ coefficient}
\\& &\quad {\rm associated\ to}\ x\  {\rm is\ equal}\ {\rm either\ to}\ \infty
\ {\rm or\ to\ a\ pair\ of\ coprime\ integers}\big\}.
\end{eqnarray*}

\begin{teo}\label{Dehn:fill:eq:teo}
If $U$ is sufficiently small and $x$ belongs
to $I\Omega_{g,k}\cap U$, then
$\widehat{M}(x)$ admits a complete finite-volume smooth
hyperbolic structure which is obtained by adding to $M(x)$ a closed
geodesic at any cusp with non-infinite Dehn filling coefficient.
From a topological point of view,
$\widehat{M} (x)$ is obtained by filling the $i$-th cusp
of $M$ along the slope $p_i (x) \mu_i+q_i (x)\lambda_i$ if
$(p_i (x),q_i (x))\neq\infty$, and by leaving the $i$-th cusp
of $M$ unfilled if $(p_i (x),q_i (x))=\infty$, $i=1,\ldots,k$.
\end{teo}
\dimostraz
See \emph{e.~g.}~\cite{Thu:bibbia,NeuZag,BenPet:book,Fri:tesi}.
\finedimo

\begin{prop}\label{added:geo:prop}
Let $U$ be sufficiently small, take $x\in I\Omega_{g,k}\cap U$ and suppose
$(p_j(x),q_j (x))\neq \infty$.
Let $\ell_j\subset \widehat{M}(x)$ be the added geodesic at the $j$-th cusp of 
$M$ and let $\matC L^x (\ell_j)$ be its complex length.
Choose integers $r_j (x), s_j (x)$ with $p_j(x)s_j (x)-q_j (x) r_j (x)=-1$.
Then we have
\[
\matC L^x (\ell_j)=r_j (x) u_j(x)+s_j (x) v_j(x).
\]
\end{prop}
\dimostraz See~\cite{NeuZag}.

\section{Isolation of cusps}\label{isolation:sec}
We now study small
deformations of the complete hyperbolic structure of $M$
by analyzing deformations of the shapes of the
geometric tetrahedra of $\calT$. 
For $x\in\Omega_{g,k}$ let $\ell(x)\in\matR$
be the length of any boundary edge of any compact lateral hexagon
in the geometric realization of $\calT$ parameterized by $x$. 
The following results are proved in~\cite[Section 3]{Fri:preprint}.

\begin{lemma}\label{compact:deform2:cor}
The map  $\ell:\Omega_{g,k}\to\matR$ is smooth, and 
$d\ell_{x_0}=0$.
\end{lemma}

\begin{lemma}\label{compact:deform:cor}
For $x\in\Omega_{g,k}$ let
$\beta (x)=x^{12k+1}$ be the 
dihedral angle along the edges of 
any compact tetrahedron in the geometric realization
of $\calT$ parameterized by $x$. 
Then 
$d\beta_{x_0}=0$.
\end{lemma}

\subsection{Infinitesimal deformations}
We begin by looking for
explicit equations for the tangent
space $T_{x_0} \Omega_{g,k}$. 
So fix a smooth arc $\varphi:(-\varepsilon,\varepsilon)\to\Omega_{g,k}$
and for any smooth $f:\Omega_{g,k}\to\matR$ let us denote
by $\dot{f}$ the derivative of $f\circ\varphi$ at $t=0$.
With notation as 
in Subsection~\ref{hyp:subsection}, 
if $\theta(t)$ is the geometric realization of $\calT$ 
parameterized by $\varphi(t)\in\Omega_{g,k}$, we set
$\alpha_l^j (t)=\theta(t)(e_l^j)$,
$\gamma_l^j(t)=\theta(t)(f_l^j)$. 
 
Recall that for every $l=1,\ldots,2k$, $j=1,2,3$ we have 
\[
\ell(\varphi(t)) =\frac{\cos\alpha_l^{j+1}(t)\cdot
\cos\alpha_l^{j+2}(t)+
\cos\gamma_l^j (t)}{\sin\alpha_l^{j+1}(t)\cdot \sin\alpha_l^{j+2}(t)},
\]
where apices are considered ${\rm mod}\ 3$. 
Moreover, by 
Lemma~\ref{compact:deform2:cor}
we have $\dot{\ell}=0$,
so differentiating at $0$ the equation above we easily get
\[
\sqrt{3} (\dot{\alpha}_l^{j+1}+
\dot{\alpha}_l^{j+2})\cos\overline{\alpha}_{g,k}+
\dot{\gamma}_l^j\sin\overline{\alpha}_{g,k}=0.
\]
Summing up these equations for $j=1,2,3$ and observing that
$\dot{\gamma}_l^1+\dot{\gamma}_l^2+\dot{\gamma}_l^3=0$ we obtain
\begin{equation}\label{euna:eq}
\dot{\alpha}_l^1+\dot{\alpha}_l^2+\dot{\alpha}_l^3=0,
\end{equation}
whence
\begin{equation}\label{edue:eq}
\sqrt{3}\dot{\alpha}_l^j \cos\overline{\alpha}_{g,k}=
\dot{\gamma}_l^j \sin\overline{\alpha}_{g,k}.
\end{equation}

Let $i\in\{1,\ldots,k\}$.
Evaluating equations~(\ref{sigma:eq}) along $\varphi$ and differentiating at $0$ we get
\[
\begin{array}{lll}
& &
\sqrt{3}(\dot{\alpha}_{2i-1}^1+\dot{\alpha}_{2i}^1)\cos\overline{\alpha}_{g,k}+
(\dot{\gamma}_{2i-1}^1+\dot{\gamma}_{2i}^1)\sin\overline{\alpha}_{g,k}\\ &=&
\sqrt{3}(\dot{\alpha}_{2i-1}^2+\dot{\alpha}_{2i}^2)\cos\overline{\alpha}_{g,k}+
(\dot{\gamma}_{2i-1}^2+\dot{\gamma}_{2i}^2)\sin\overline{\alpha}_{g,k}\\ &=&
\sqrt{3}(\dot{\alpha}_{2i-1}^3+\dot{\alpha}_{2i}^3)\cos\overline{\alpha}_{g,k}+
(\dot{\gamma}_{2i-1}^3+\dot{\gamma}_{2i}^3)\sin\overline{\alpha}_{g,k}.
\end{array}
\]
Together with equations~(\ref{euna:eq}) 
and~(\ref{edue:eq}), this implies 
\begin{equation}\label{etre:eq}
\dot{\alpha}_{2i-1}^1=-\dot{\alpha}_{2i}^1,\quad
\dot{\alpha}_{2i-1}^2=-\dot{\alpha}_{2i}^2,\quad
\dot{\alpha}_{2i-1}^3=-\dot{\alpha}_{2i}^3.
\end{equation}

We can now summarize all these computations giving explicit
equations for $T_{x_0} \Omega_{g,k}$. 
Let $Z$ be the linear subspace of $\matR^{12}$
defined by the following equations: 
\[
\left\{
\begin{array}{l}
(\sqrt{3}\cos\overline{\alpha}_{g,k})x_1=
(\sin\overline{\alpha}_{g,k})x_4\\
(\sqrt{3}\cos\overline{\alpha}_{g,k})x_2=
(\sin\overline{\alpha}_{g,k})x_5\\
(\sqrt{3}\cos\overline{\alpha}_{g,k})x_3=
(\sin\overline{\alpha}_{g,k})x_6\\
x_1+x_7=x_2+x_8=x_3+x_9=x_4+x_{10}=x_5+x_{11}=x_6+x_{12}=0\\
x_1+x_2+x_3=0.
\end{array}
\right.
\]
Observe that $\dim_\matR Z=2$.
For $i=1,\ldots, k$ let $r_i:\matR^{12k+1}\to\matR^{12}$
be the map defined by $r_i(x)=(x_{12i-11},x_{12i-10},\ldots,
x_{12i-1},x_{12i})$. Let $$\overline{Z}=\{x\in\matR^{12k+1}:\ 
x_{12k+1}=0,\, r_i(x)\in Z\ {\rm for}\ i=1,\ldots,k\}$$ 
be the product of one copy of 
$Z$ for each cusp (so $\dim_\matR Z=2k$).

\begin{prop}\label{last:step:prop}
We have $T_{x_0}\Omega_{g,k}=\overline{Z}$.
\end{prop}
\dimostraz
Lemma~\ref{compact:deform:cor} 
and equations~(\ref{euna:eq}),~(\ref{edue:eq}),~(\ref{etre:eq}),
imply that $T_{x_0}\Omega_{g,k}\subseteq \overline{Z}$.
But $\dim_\matR \overline{Z}= 2k=\dim_\matR T_{x_0}\Omega_{g,k}$, whence the
conclusion.
\finedimo

\subsection{Isolation of cusps}
We now go into the proof of Theorem~\ref{isolation:teo}.
So let $C_1,\ldots,C_k$ be the cusps of our fixed manifold
$M\in\calM_{g,k}$ corresponding to the boundary tori $T_1,\ldots,T_k$.
We look for
equations defining the set of structures in $\Omega_{g,k}$
which are complete at $C_1,\ldots,C_h$. 
To this aim we set:
\begin{equation}\label{defh:eq}
J_h=\{x\in\matR^{12k+1}:\ x_{12i+1}=x_{12i+2}=x_{12i+3}
\ {\rm for\ all}\ i=0,\ldots,h-1\}.
\end{equation}

\begin{lemma}\label{iso1:lemma}
Near $x_0$, the set  $J_h\cap \Omega_{g,k}$
is a smooth submanifold of $\Omega_{g,k}$ 
of real dimension $2(k-h)$.
\end{lemma}
\dimostraz
It is easily seen that $T_{x_0} J_h +T_{x_0} \Omega_{g,k}=
J_h +\overline{Z}=\matR^{12k+1}=T_{x_0} \matR^{12k+1}$, so 
the conclusion follows from
basic results about transverse intersections
of submanifolds.
\finedimo

Let $\Delta$ be a topological
partially truncated tetrahedron with ideal
vertex $v_0$, and take $\vartheta\in (0,\pi/3)$. Then
there exists, up to isometry, exactly one geometric
realization of $\Delta$ with dihedral angles $\pi/3$
along the internal edges emanating
from $v_0$, and angle $\vartheta$ along the other internal edges.
We denote this geometric tetrahedron by $\Delta^\vartheta$.

\begin{prop}\label{iso2:prop}
For each $l=1,\ldots,2h$
let $\Delta^{\!\ast}_l(p)$ be the geometric realization
of $\Delta_l$ parameterized by $p\in J_h \cap\Omega_{g,k}$.
Then a real number $\vartheta (p)\in (0,\pi/3)$ exists
such that $\Delta^{\!\ast}_l(p)$ is isometric to
$\Delta^{\vartheta (p)}$ for $l=1,\ldots,2h$.
\end{prop}
\dimostraz
For $l=1,\ldots,2h$, $j=1,2,3$  
let $T^j_l(p)$ be the truncation triangle
of $\Delta^{\!\ast}_l(p)$ having a vertex on
the edge $f_l^j$.
Fix $i\in\{0,\ldots,h-1\}$ and consider the tetrahedron
$\Delta^{\!\ast}_{2i+1} (p)$. The compact face of such tetrahedron
is a \emph{regular} right-angled hexagon,
so condition $x_{12i+1} (p)=x_{12i+2}(p)=x_{12i+3}(p)$ implies that
$T^1_{2i+1}(p),T^2_{2i+1}(p)$ and $T^3_{2i+1}(p)$ 
are isometric to each other. This 
gives $x_{12i+4} (p)=x_{12i+5}(p)=x_{12i+6}(p)$, so
$\Delta^{\!\ast}_{2i+1}$ is isometric to $\Delta^{\xi_{2i+1}}$ for some
$\xi_{2i+1}=x_{12i+1} (p)\in (0,\pi/3)$. Moreover, since the non-compact faces
of $\Delta^{\!\ast}_{2i+1}$ are isometrically glued to
the non-compact faces of $\Delta^{\!\ast}_{2i+2}$ we easily see
that the truncation triangles 
$T^1_{2i+2}(p),T^2_{2i+2}(p)$ and $T^3_{2i+2}(p)$ 
are isosceles and isometric to each other. This forces
$x_{12i+7} (p)=x_{12i+8}(p)=x_{12i+9}(p)=\xi_{2i+2}$ and
$x_{12i+10} (p)=x_{12i+11}(p)=x_{12i+12}(p)=\pi/3$, so
$\Delta^{\!\ast}_{2i+2}$ is isometric to
$\Delta^{\xi_{2i+2}}$ for some
$\xi_{2i+2}\in (0,\pi/3)$.
  
Finally, since
the length
of the compact internal edges of the $\Delta^{\!\ast}_l$'s
does not depend on $l$, we have $\xi_{1}=\ldots=\xi_{2h}$,
whence the conclusion.
\finedimo

\begin{cor}\label{iso3:cor}
Let $p$ be a point in $J_h \cap\Omega_{g,k}$ and denote
by $M(p)$ the hyperbolic structure defined by $p$ on $M$.
Then for all $i=1,\ldots,h$ the following holds:
\begin{itemize}
\item
$M(p)$ induces a complete metric on the cusp $C_i$;
\item
The Euclidean structure induced on $T_i$
by $M(p)$ is isometric to the regular hexagonal structure
induced on $T_i$ by the complete hyperbolic structure
$M(x_0)$.
\end{itemize}
\end{cor}

The corollary just stated says that the Euclidean structures
on $T_1,\ldots, T_h$ are not affected by the deformations of the
hyperbolic metric on $M$ which correspond to points
in $J_h \cap \Omega_{g,k}$. Therefore to conclude the proof of
Theorem~\ref{isolation:teo} we only need the following:

\begin{prop}\label{iso4:prop}
Let $K_h $ be the subset of $\Omega_{g,k}$ corresponding
to the structures inducing \emph{complete} metrics on
$C_1,\ldots,C_h$. Then there exists a neighbourhood $V$
of $x_0$ in $\Omega_{g,k}$ with
$K_h \cap V=J_h \cap V$.
\end{prop}
\dimostraz
By Theorem~\ref{smoothness:teo}
and Lemma~\ref{iso1:lemma} there exists a  neighbourhood 
$W$ of $x_0$ in $\Omega_{g,k}$ such that both $K_h \cap W$ and $J_h \cap W$
are smooth submanifolds of $\Omega_{g,k}$ of real dimension $2(k-h)$. 
Moreover Corollary~\ref{iso3:cor} shows that $J_h \cap W\subset K_h \cap W$,
whence $J_h \cap V=K_h \cap V$ for some (maybe smaller) neighbourhood
of $x_0$ in $\Omega_{g,k}$.
\finedimo

\begin{rem}
Theorem~\ref{isolation:teo} shows that 
Dehn fillings along \emph{sufficiently complicated}
slopes on some boundary tori of $M$ 
do not affect the Euclidean structure on the non-filled 
boundary tori. Using SnapPea we have checked in a number of cases
that the same isolation phenomenon still holds when filling along
\emph{short}
non-exceptional slopes. It is   
conjectured in~\cite{HodKer} that the space of Dehn filling
deformations of complete finite-volume
hyperbolic $3$-manifolds is connected
and smooth.
If this were true, then the Euclidean structure on the non-filled 
tori would remain unchanged under all the partial hyperbolic
Dehn fillings of $M$.
\end{rem}

\section{Non-isolation of the boundary}\label{nonisolation:section}

Let ${\rm Teich} (\partial M)$ be
the Teichm\"uller space of hyperbolic structures on $\partial M$,
\emph{i.~e.~} the space of equivalence classes of hyperbolic metrics
on $\partial M$, where two such metrics are considered equivalent if
they are isometric through a diffeomorphism homotopic to the identity
of $\partial M$. Since $\partial M$ is compact, it is well-known that
for any $\gamma\in\pi_1 (\partial M)$ and any metric $h\in {\rm Teich} (\partial M)$
there exists a unique closed
$h$-geodesic in the free homotopy class of $\gamma$.
We denote the $h$-length of this geodesic by $L_\gamma (h)$. An easy computation
shows that 
if $\overline{\rho}_h:\pi_1 (\partial M)\to {\rm PSL}(2,\matR)\cong {\rm Isom}^+ (\matH^2)$
is a holonomy representation for the hyperbolic structure $h$ then we
have
\begin{equation}\label{trace:length:eq}
{\rm tr}\, \rho_h (\gamma)=\pm 2 \cosh (L_\gamma (h)/2),
\quad \gamma\in\pi_1 (\partial M), 
\end{equation}
where $\rho_h (\gamma)$ is a lift in ${\rm SL}(2,\matR)$ of 
$\overline{\rho}_h (\gamma)\in {\rm PSL}(2,\matR)$. 

For $x\in\Omega_{g,k}$ we denote by $M(x)$ the hyperbolic
structure defined on $M$ by $x$, and by $B(x)\in {\rm Teich}
(\partial M)$ the equivalence class of the hyperbolic structure
induced by $M(x)$ on $\partial M$.
It is well-known that ${\rm Teich} (\partial M)$
admits a structure of differentiable manifold such that:
\begin{itemize}
\item
${\rm Teich} (\partial M)$ is diffeomorphic to the Euclidean
space
$\matR^{6g-6}$;
\item
For any $\gamma\in\pi_1 (M)$ the map 
$L_\gamma:{\rm Teich} (\partial M)\to\matR$
defined above
is smooth;
\item
The map $B:\Omega_{g,k}\to {\rm Teich} (\partial M)$
is smooth. 
\end{itemize}

The following proposition is proved in~\cite{Fri:preprint}:

\begin{prop}\label{diff0:prop}
We have $dB_{x_0}=0$.
\end{prop} 

Thus in order to understand how deformations of the complete structure
affect the geodesic boundary we need to analyze the map
$B$ up to the second order in a neighbourhood of $x_0$.
We begin with the following:

\begin{defn}\label{direction:defn}
Let $y_0$ be a point of a smooth $n$-manifold $Y$ and let
$\varphi:U\to\matR^n$ be a diffeomorphism with $\varphi (y_0)=0$, where
$U\subset Y$ is a small open neighbourhood of $y_0$.
Let $0\neq v\in T_{y_0} Y$, and consider
a sequence $\{y_j\}_{j\in\matN}\subset U\setminus\{y_0\}$. We say that $y_n$
converges to $y_0$ \emph{along $v$} if
\[
\lim_{j\to\infty} y_j= y_0, \quad
\lim_{j\to\infty} \varphi(y_j)/||\varphi(y_j)||= d\varphi_{y_0} (v)/
||d\varphi_{y_0} (v)||,
\]
where we are identifying $T_0 \matR^n$ with $\matR^n$, endowed
with  the Euclidean norm $||\cdot ||$.
\end{defn}

\subsection{An alternative formulation}
First of all we show how Theorem~\ref{nonisol:intro:teo}
can be deduced from the 
following:

\begin{teo}\label{nonisol:key:teo}
There exist 
a smooth path $\overline{\varsigma}:
(-\varepsilon,\varepsilon)\to \Omega_{g,k}$ 
and an element $\overline{\gamma}\in\pi_1 (\partial M)$ 
such that
$\overline{\varsigma} (0)=x_0$ and
the map $t\mapsto L^{B (\overline{\varsigma} (t))} (\overline{\gamma})$ has
non-zero second derivative at $0$.
\end{teo}

Let $\{y_n\}_{n\in\matN}\subset  I\Omega_{g,k}\setminus
\{x_0\}$ be a sequence converging to $x_0$ along
$\dot{\overline{\varsigma}}(0)$
(the set 
$I\Omega_{g,k}\subset \Omega_{g,k}$ was defined in 
Subsection~\ref{Dehn:subsection}).
By construction, up to passing to a subsequence
we have $B(y_n)\neq B(x_0)$ for every $n\in\matN$.

Since $y_i\in I\Omega_{g,k}$,
the metric completion of the structure
induced on $M$ by $y_i$ gives a non-singular hyperbolic
$3$-manifold $N_i$.
Recall that the mapping class group
$\mathcal{MCG} (\partial M)$  of $\partial M$
acts properly discontinuously on
${\rm Teich} (\partial M)$, so there exists 
a neighbourhood $U$
of $B(x_0)$ in ${\rm Teich} (\partial M)$ such that
the set $\{\psi\in\mathcal{MCG} (\partial M):\,
\psi (U)\cap U\neq \emptyset\}$ is finite.
Up to passing to a subsequence,
we may suppose that the equivalence classes of the $\partial N_j$'s
are pairwise distinct as elements in ${\rm Teich} (\partial M)$, and that 
$\partial N_i\in U$ for all $i\in\matN$. This readily implies that
among the $\partial N_j$'s there are infinitely many pairwise non-isometric
hyperbolic surfaces, whence Theorem~\ref{nonisol:intro:teo}.

\subsection{Proving Theorem~\ref{nonisol:key:teo}}
Let $\varsigma:(-\varepsilon,\varepsilon)\to\Omega_{g,k}$ 
be a fixed smooth path with $\varsigma(0)=x_0$ and for any smooth function
$f:\Omega_{g,k}\to\matR$ let $\dot{f}(t)$ (resp.~$\ddot{f}(t)$)
denote the first (resp.~second) derivative of $f\circ\varsigma$ in $t$.
We will denote by $\dot{f}$
(resp.~$\ddot{f}$) 
the value of $\dot{f}(0)$ (resp.~of $\ddot{f}(0)$).
We recall that $\dot{\ell}=\dot{x}_{{12k+1}}=0$
by Lemmas~\ref{compact:deform2:cor},~\ref{compact:deform:cor}.

\begin{lemma}\label{second:length:lemma}
We have 
$\ddot{\ell}=\ddot{x}_{6l+1}+\ddot{x}_{6l+2}+
\ddot{x}_{6l+3}=\ddot{x}_{{12k+1}}=0$ for
$l=0,\ldots,2k-1$.
\end{lemma}
\dimostraz
Differentiating two times the equality
$\ell(x)=(\cos x_1\cos x_2+\cos x_6)/(\sin x_1 \sin x_2)$
along $\varsigma$ 
and evaluating at $0$ we get 
\begin{equation}\label{dersec1:eq}
\begin{array}{llll}
(2\sin^4 \overline{\alpha}_{g,k})\ddot{\ell}&=&
(-3\cos\overline{\alpha}_{g,k}\sin\overline{\alpha}_{g,k})(\ddot{x}_{1}+
\ddot{x}_{2})-\sqrt{3}\sin^2 \overline{\alpha}_{g,k}\ddot{x}_{6}\\ &+&
(5\cos^2 \overline{\alpha}_{g,k}+1) (\dot{x}_{1}^2+\dot{x}_{2}^2)
-\sin^2 \overline{\alpha}_{g,k} \dot{x}_{6}^2\\ &+&
(2\cos^2 \overline{\alpha}_{g,k}+4)\dot{x}_{1}\dot{x}_{2}\\ &+&
2\sqrt{3}\cos \overline{\alpha}_{g,k} \sin \overline{\alpha}_{g,k}
(\dot{x}_{1}\dot{x}_{6}+\dot{x}_{2}\dot{x}_{6}).
\end{array}
\end{equation}
Observe now that on $\Omega_{g,k}$ we also have $$\ell (x)=
(\cos x_1\cos x_3+\cos x_5)/
(\sin x_1\sin x_3)=(\cos x_2\cos x_3+\cos x_4)/
(\sin x_2\sin x_3).$$ 
Differentiating two times these equalities along $\varsigma$,
evaluating at $0$ as above 
and 
summing up with equality~(\ref{dersec1:eq}) we get
\begin{equation}\label{finald2:eq}
\begin{array}{llll}
(6\sin^4 \overline{\alpha}_{g,k})\ddot{\ell}&=&  
-6\cos\overline{\alpha}_{g,k}\sin\overline{\alpha}_{g,k}
(\ddot{x}_{1}+\ddot{x}_{2}+\ddot{x}_{3})-\sqrt{3}\sin^2
\overline{\alpha}_{g,k}(\ddot{x}_{4}+\ddot{x}_{5}+\ddot{x}_{6})\\
&+& (10\cos^2 \overline{\alpha}_{g,k}+2)(\dot{x}_{1}^2
+\dot{x}_{2}^2+\dot{x}_{3}^2)-\sin^2 \overline{\alpha}_{g,k}
(\dot{x}_{4}^2+\dot{x}_{5}^2+\dot{x}_{6}^2)\\
&+& (2\cos^2\overline{\alpha}_{g,k} +4)(\dot{x}_{1}\dot{x}_{2}
+\dot{x}_{1}\dot{x}_{3}+\dot{x}_{2}\dot{x}_{3})\\&+&
2\sqrt{3}\cos\overline{\alpha}_{g,k}\sin\overline{\alpha}_{g,k}
(\dot{x}_{1}\dot{x}_{5}+\dot{x}_{1}\dot{x}_{6}+\dot{x}_{2}
\dot{x}_{4}+\dot{x}_{2}\dot{x}_{6}+\dot{x}_{3}\dot{x}_{4}
+\dot{x}_{3}\dot{x}_{5}).
\end{array}
\end{equation}

Since $\dot{x}$ belongs to $T_{x_0} \Omega_{g,k}$ we have  
$(\sin\overline{\alpha}_{g,k})\dot{x}_{{i+3}}=
(\sqrt{3}\cos\overline{\alpha}_{g,k})\dot{x}_{i}$ for $i=1,2,3$ and
$\dot{x}_{1}+\dot{x}_{2}+\dot{x}_{3}=0$. 
Substituting these relations
in~(\ref{finald2:eq}) we (rather strikingly) get
$\ddot{\ell}=-(\cos\overline{\alpha}_{g,k}/\sin^3 \overline{\alpha}_{g,k})
(\ddot{x}_{1}+
\ddot{x}_{2}+\ddot{x}_{3})$. 
By the very same argument
it follows that 
\begin{equation}\label{third:d:eq}
\ddot{\ell}=-(\cos\overline{\alpha}_{g,k}/\sin^3 \overline{\alpha}_{g,k})
(\ddot{x}_{6l+1}+
\ddot{x}_{6l+2}+\ddot{x}_{6l+3}), 
\quad l=0,\ldots,2k-1. 
\end{equation}
Thus the condition forcing the dihedral angle along the 
compact edge of $\calT$ to
be constantly equal to $2\pi$ gives
$-2k \ddot{\ell} \sin^3\overline{\alpha}_{g,k}/\cos\overline{\alpha}_{g,k}
+6(g-k)\ddot{x}_{{12k+1}}=0$, and implies that
$\ddot{\ell}$ has the same sign as
$\ddot{x}_{{12k+1}}$.
On the other hand, condition
$\ell(x)=\cos x_{12k+1} /
(1-\cos x_{12k+1})$ implies $\ddot{\ell}=
-(\sin\overline{\beta}_{g,k}/(1-\cos\overline{\beta}_{g,k})^2)
\ddot{x}_{{12k+1}}$, so 
$\ddot{\ell}$
and $\ddot{x}_{{12k+1}}$ should have opposite signs.
This forces $\ddot{\ell}=\ddot{x}_{12k+1}=0$, whence the conclusion
by equation~(\ref{third:d:eq}).
\finedimo

\subsection{The chosen curve}\label{curve:subsection}
By Proposition~\ref{last:step:prop}
the subspace of $\matR^{12k+1}$ having equations
$\{x\in\matR^{12k+1}:\, x_2=x_3,\, x_{12i+1}=x_{12i+2}=x_{12i+3},
\, i=1,\ldots,k-1\}$
intersects
$\Omega_{g,k}$ transversely near $x_0$ in the support of a
smooth curve $\overline{\varsigma}:(-\epsilon,\epsilon)\to \Omega_{g,k}$ with
$\overline{\varsigma}(0)=x_0$ and
\begin{equation}\label{der:sig:eq}
\begin{array}{lll}
\dot{\overline{\varsigma}}(0) &=&
(2\sin\overline{\alpha}_{g,k},-\sin\overline{\alpha}_{g,k},
-\sin\overline{\alpha}_{g,k},2\sqrt{3}\cos\overline{\alpha}_{g,k},
-\sqrt{3}\cos\overline{\alpha}_{g,k},-\sqrt{3}\cos\overline{\alpha}_{g,k},\\
& &-2\sin\overline{\alpha}_{g,k},\sin\overline{\alpha}_{g,k},
\sin\overline{\alpha}_{g,k},-2\sqrt{3}\cos\overline{\alpha}_{g,k},
\sqrt{3}\cos\overline{\alpha}_{g,k},\sqrt{3}\cos\overline{\alpha}_{g,k},\\
& &0,\ldots,0).
\end{array}
\end{equation}
As before, for any smooth $f:\Omega_{g,k}\to\matR$ we set
$f(t):=f(\overline{\varsigma}(t))$ and we denote by 
$\dot{f}(t)$ (resp.~$\ddot{f}(t)$) the first (resp.~second)
derivative of $f\circ\overline{\varsigma}$ in $t$. We
also denote simply by $\dot{f}$, $\ddot{f}$
the values $\dot{f}(0)$, $\ddot{f}(0)$ respectively.

If $\psi:(-\epsilon',\epsilon')\to(-\epsilon,\epsilon)$ is a 
local diffeomorphism with $\dot{\psi}(0)=1$, then 
$\ddot{(\overline{\varsigma}\circ\psi)}(0)=\ddot{\overline{\varsigma}}(0)
+\ddot{\psi}(0)\cdot \dot{\overline{\varsigma}} (0)$.
Thus, up to reparameterizing $\overline{\varsigma}$ without changing
its tangent vector at $0$ we can assume that the following 
condition holds:
\begin{equation}\label{cond:der:eq}
\ddot{x}_1=\ddot{x}_7.
\end{equation} 

A very similar argument to the proof of
Proposition~\ref{iso2:prop} gives the following:

\begin{prop}\label{isol:first:prop}
For any $t\in(-\epsilon,\epsilon)$ we have
\[
x_2 (t)=x_3 (t),\quad x_5 (t)=x_6 (t),
\quad x_8(t)=x_9(t),\quad x_{11}(t)=x_{12}(t),
\]
and
\[
\begin{array}{ll}
x_{12i+1}(t)=x_{12i+2}(t)=x_{12i+3}(t),&
x_{12i+4}(t)=x_{12i+5}(t)=x_{12i+6}(t)=\pi/3,\\
x_{12i+7}(t)=x_{12i+8}(t)=x_{12i+9}(t),&
x_{12i+10}(t)=x_{12i+11}(t)=x_{12i+12}(t)=\pi/3
\end{array}
\] 
for $i=1,\ldots,k-1$.
\end{prop}

We are now ready to prove the following:

\begin{prop}\label{fine:der:prop}
We have 
\[
\begin{array}{lll}
\ddot{\overline{\varsigma}}(0)&=&
(8\cos\overline{\alpha}_{g,k}\sin\overline{\alpha}_{g,k},
-4\cos\overline{\alpha}_{g,k}\sin\overline{\alpha}_{g,k},
-4\cos\overline{\alpha}_{g,k}\sin\overline{\alpha}_{g,k},
2\sqrt{3},-\sqrt{3},-\sqrt{3},\\
& &
8\cos\overline{\alpha}_{g,k}\sin\overline{\alpha}_{g,k},
-4\cos\overline{\alpha}_{g,k}\sin\overline{\alpha}_{g,k},
-4\cos\overline{\alpha}_{g,k}\sin\overline{\alpha}_{g,k},
2\sqrt{3},-\sqrt{3},-\sqrt{3},\\
& & 0,\ldots,0).
\end{array}
\]
\end{prop}
\dimostraz
Since $x_4(t)+x_5(t)+x_6(t)=x_{10}(t)+x_{11}(t)+x_{12}(t)=\pi$,
by Proposition~\ref{isol:first:prop} 
we have $\ddot{x}_5=\ddot{x}_6=-(1/2)\ddot{x}_4$ and
$\ddot{x}_{11}=\ddot{x}_{12}=-(1/2)\ddot{x}_{10}$,
while by Lemma~\ref{second:length:lemma} 
and Proposition~\ref{isol:first:prop} we get
$\ddot{x}_2=\ddot{x}_3=-(1/2)\ddot{x}_1$ and
$\ddot{x}_{8}=\ddot{x}_{9}=-(1/2)\ddot{x}_{7}$.
Substituting these relations in equation~(\ref{dersec1:eq})
and recalling that $\ddot{\ell}=0$ 
we have after some computations
\begin{equation}\label{derd1:eq}
\sqrt{3}\cos\overline{\alpha}_{g,k}\ddot{x}_1
-\sin\overline{\alpha}_{g,k}\ddot{x}_4=
2\sqrt{3}\sin\overline{\alpha}_{g,k}
(4\cos^2\overline{\alpha}_{g,k}-1).
\end{equation}
The very same argument also gives
\begin{equation}\label{derd2:eq}
\sqrt{3}\cos\overline{\alpha}_{g,k}\ddot{x}_7
-\sin\overline{\alpha}_{g,k}\ddot{x}_{10}=
2\sqrt{3}\sin\overline{\alpha}_{g,k}
(4\cos^2\overline{\alpha}_{g,k}-1).
\end{equation}
In the same way, differentiating two
times the equality 
\[
\sin x_{1}(t)\sin x_{4}(t)\sin x_{7}(t)\sin x_{10}(t)=
\sin x_{2}(t)\sin x_{5}(t)\sin x_{8}(t)\sin x_{11}(t),
\]
evaluating at $0$ and substituting in the result 
the relations $\ddot{x}_2=-(1/2)\ddot{x}_1$, 
$\ddot{x}_5=-(1/2)\ddot{x}_4$,
$\ddot{x}_{8}=-(1/2)\ddot{x}_{7}$ and
$\ddot{x}_{11}=-(1/2)\ddot{x}_{10}$
we get
\begin{equation}\label{derd3:eq}
\sqrt{3}\cos\overline{\alpha}_{g,k}(\ddot{x}_1+\ddot{x}_7)
+\sin\overline{\alpha}_{g,k}(\ddot{x}_4+\ddot{x}_{10})=
4\sqrt{3}\sin\overline{\alpha}_{g,k}(4\cos^2\overline{\alpha}_{g,k}+1).
\end{equation}
Solving equations~(\ref{cond:der:eq}), (\ref{derd1:eq}), (\ref{derd2:eq})
and
(\ref{derd3:eq}) we get the desired
result for $\ddot{x}_1,\ldots,\ddot{x}_{12}$.

Let now $i=1,\ldots,k-1$. By Proposition~\ref{isol:first:prop}
and
Lemma~\ref{second:length:lemma}
we get 
\begin{eqnarray*}
\ddot{x}_{12i+1}=\ddot{x}_{12i+2}=\ddot{x}_{12i+3}=
(\ddot{x}_{12i+1}+\ddot{x}_{12i+2}+\ddot{x}_{12i+3})/3=0,\\
\ddot{x}_{12i+7}=\ddot{x}_{12i+8}=\ddot{x}_{12i+9}=
(\ddot{x}_{12i+7}+\ddot{x}_{12i+8}+\ddot{x}_{12i+9})/3=0.
\end{eqnarray*}
Moreover, Proposition~\ref{isol:first:prop} forces
\[
\ddot{x}_{12i+4}=\ddot{x}_{12i+5}=\ddot{x}_{12i+6}=\ddot{x}_{12i+10}=
\ddot{x}_{12i+10}=\ddot{x}_{12i+11}=\ddot{x}_{12i+12}=0,
\]
and 
Lemma~\ref{second:length:lemma} also gives
$\ddot{x}_{12k+1}=0$. 
\finedimo

\begin{rem}\label{new:ext:rem}
For $i=1,\ldots, k$ recall that the map $r_i:\matR^{12k+1}\to\matR^{12}$
is defined by $r_i(x)=(x_{12i-11},x_{12i-10},\ldots,
x_{12i-1},x_{12i})$. 
If $h\leqslant k$, a slight modification of
the strategy adopted to construct
$\overline{\varsigma}$ yields
a curve $\overline{\varsigma}_h:(-\varepsilon,\varepsilon)\to\Omega_{g,k}$
with the following properties: 
$\overline{\varsigma}_h (0)=x_0$, 
the structure defined by $\overline{\varsigma} (t)$ on $M$ 
is complete at the last $k-h$ cusps, and
$r_i (\dot{\overline{\varsigma}}_h(0))
=r_1 (\dot{\overline{\varsigma}} (0))$,  
$r_i (\ddot{\overline{\varsigma}}_h(0))
=r_1 (\ddot{\overline{\varsigma}} (0))$ for all $i=1,\ldots,h$.
\end{rem}

\subsection{The final step}
The smooth path $\overline{\varsigma}:(-\varepsilon,\varepsilon)\to\Omega_{g,k}$ 
determines a smooth family of developing maps
$D_t: \widetilde{\partial M}\to\matH^2$, which gives in turn
a smooth path of holonomy representations
$\overline{\rho}_t:\pi_1 (\partial M)\to {\rm PSL}(2,\matR)$ lifting
to a smooth path of representations $\rho_t:\pi_1 (\partial M)
\to {\rm SL}(2,\matR)$. For any $\gamma\in \pi_1 (\partial M),\, t\in
(-\varepsilon,\varepsilon)$ we set
${\rm tr}_\gamma (t)={\rm trace}\,\rho_t (\gamma)$.
Of course ${\rm tr}_\gamma:(-\varepsilon,\varepsilon) \to \matR$
is smooth for any $\gamma\in\pi_1 (\partial M)$.
Moreover, 
from Proposition~\ref{diff0:prop}
and equation~(\ref{trace:length:eq}) 
we easily deduce $\dot{\tr}_\gamma (0)=0$.
The following result readily implies 
Theorem~\ref{nonisol:key:teo}, whence
Theorem~\ref{nonisol:intro:teo}.

\begin{prop}\label{nonisol:prop}
There exists an element $\overline{\gamma}\in\pi_1 (\partial M)$
such that $\ddot{\tr}_{\overline{\gamma}} (0)\neq 0$.
\end{prop}

\begin{figure}
\begin{center}
\input{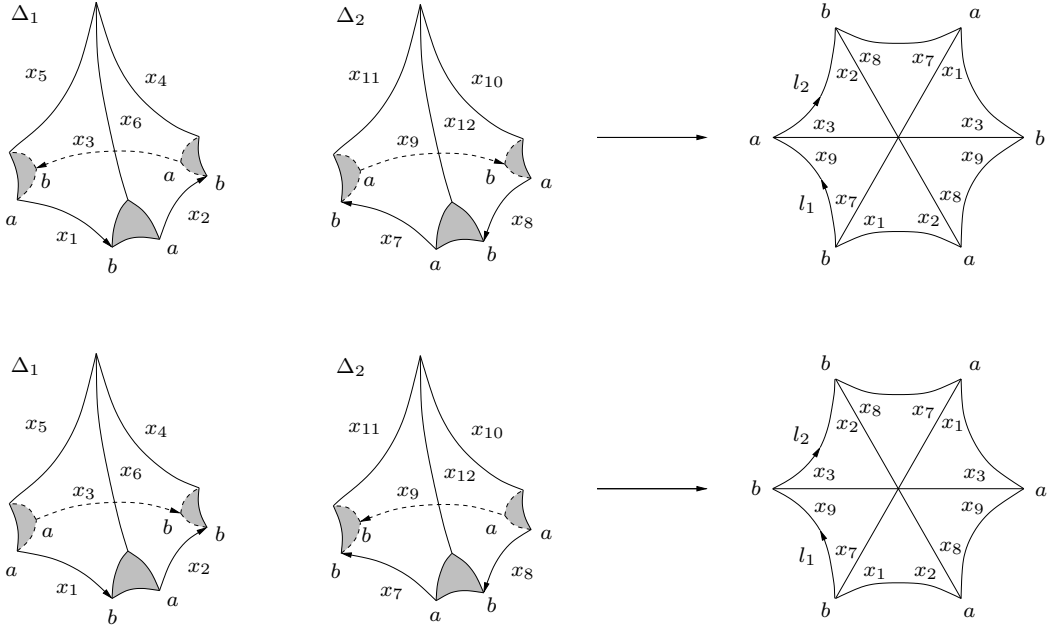}
\caption{Up to cyclic reorderings of the internal
edges $e_1^1,e_1^2,e_1^3$ of $\Delta_1$, we show here
the only 
possible identifications of the vertices of the hexagon tessellated by the
truncation triangles of $\Delta_1$ and $\Delta_2$.}\label{hexagon:fig}
\end{center}
\end{figure}

In order to find the element $\overline{\gamma}\in\pi_1 (\partial M)$
mentioned in Proposition~\ref{nonisol:prop}
we have to describe in some detail a particular portion
of the triangulation induced on $\partial M$ by the canonical
decomposition $\calT$ of $M$.
So we fix our attention on the geodesic
hexagon (with identifications on edges and vertices)
which results from the gluing of the truncation triangles of
$\Delta_1,\Delta_2\in\calT$. Let $l_1, l_2$ be the oriented
edges of this
hexagon described in Fig.~\ref{hexagon:fig}
 and observe that the starting point of $l_1$ and the endpoint of $l_2$
both coincide with the same point $b\in\partial M$.
Thus the loop $l_2 \ast l_1$ defines an element $\overline{\gamma}$
of $\pi_1 (\partial M, b)$. Our next aim is to give an explicit
description of the isometry $\overline{\rho}_t (\overline{\gamma})
\in {\rm PSL}(2,\matR)$ 
in terms of angles and lengths of edges of the triangulation
of $\partial M$. Let us fix two consecutive lifts
$\widetilde{l}_1,\widetilde{l}_2$ 
of $l_1, l_2$ in $\widetilde{\partial M}$ and let
$\widetilde{l}'_1=\overline{\gamma} (\widetilde{l}_1)\subset\widetilde{\partial M}$
be the lift of $l_1$ starting at the endpoint of $\widetilde{l}_2$.
Then $\overline{\rho}_t (\overline{\gamma})$ 
is the unique orientation-preserving isometry
of $\matH^2$ taking the oriented geodesic segment 
$D_t (\widetilde{l}_1)$ onto the oriented geodesic segment
$D_t (\widetilde{l}'_1)$.

\begin{figure}
\begin{center}
\input{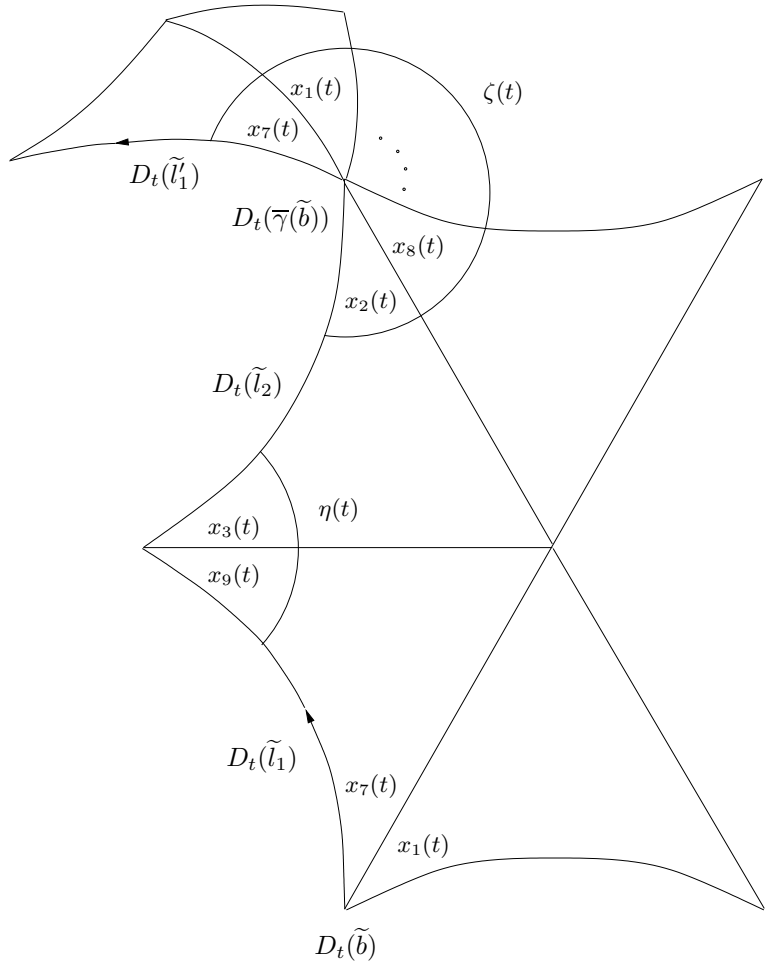}
\caption{{Notation for the proof of 
equation~(\ref{der:hol:eq}).}}\label{holbordo:fig}
\end{center}
\end{figure}

Let $\eta (t)$ be the angle formed by $D_t (\widetilde{l}_1)$ and
$D_t (\widetilde{l}_2)$ at the endpoint of $D_t (\widetilde{l}_1)$ and
$\zeta (t)$  the angle formed by $D_t (\widetilde{l}_2)$ and
$D_t (\widetilde{l}'_1)$ at the endpoint of $D_t (\widetilde{l}_2)$
(see Fig.~\ref{holbordo:fig}).
Of course we have $\eta (t)= x_3 (t)+x_9 (t)$, while
\[
\zeta (t)=(x_2 (t)+x_8 (t))+ (x_1 (t)+x_7 (t))+\delta
(x_3 (t)+x_9 (t))+ r(t),
\]
where $\delta\in\{0,1\}$ is determined by the combinatorics
of $\calT$ and $r(t)\in (0,2\pi)$ is given by the sum (with multiplicity)
of some of the $x_i (t)$'s with $i\geqslant 13$.
Equation~(\ref{der:sig:eq})
and Proposition~\ref{fine:der:prop} give the following:

\begin{lemma}\label{alphabeta:lemma}
We have $\dot{\eta}(0)=\dot{\zeta}(0)=0$ and $\ddot{\eta}(0)\neq 0$. Moreover, 
if $\delta=0$ then $\ddot{\zeta}(0)=-\ddot{\eta}(0)$;
if $\delta=1$ then $\ddot{\zeta}(0)=0$.
\end{lemma}

Let now $\matH=\{z\in\matC:\ \Im (z)>0\}$ be the upper half-plane model of
$\matH^2$. Without loss of generality we may assume that
for any $t\in(-\varepsilon,\varepsilon)$ the developing map
$D_t: \widetilde{\partial M}\to\matH$ sends $\widetilde{l}_1$
onto the geodesic segment starting at $i\in\matH$ and ending
at $\lambda(t)\cdot i\in\matH$, where $\lambda(t)=\exp \ell(t)$.
For $1<\lambda\in\matR,
\ \theta\in (0,2\pi)$ we set
\[
A(\lambda,\theta)=\left[\left(
\begin{array}{cc}
\sqrt{\lambda} \sin (\theta/2) & -\sqrt{\lambda} \cos (\theta/2)\\
(1/\sqrt{\lambda})\cos (\theta/2) & (1/\sqrt{\lambda})\sin (\theta/2)
\end{array}
\right)\right] \in {\rm PSL}(2,\matR). 
\]
It is easily seen that $A(\lambda,\theta)$ sends the half-geodesic
$s$ starting at $i\in\matH$ and ending at $\infty$ onto the half-geodesic 
$s'$ starting at $\lambda\cdot i$ such that $s$ and $s'$ define at $\lambda\cdot i$
an angle equal to $\theta$.
Now both $D_t (\widetilde{l}_2)$ 
and $D_t (\widetilde{l}_1)$ have length $\ell (t)$, and the isometry 
$\overline{\rho}_t (\overline{\gamma})$ takes the oriented 
geodesic segment $D_t (\widetilde{l}_1)$ onto the oriented 
geodesic segment $D_t (\widetilde{l}'_1)$, so 
\begin{equation}\label{holcomp:eq}
\overline{\rho}_t (\overline{\gamma})=
A(\lambda(t),\eta (t))\cdot A(\lambda(t), \zeta (t)).
\end{equation}
By
Lemmas~\ref{compact:deform2:cor},~\ref{second:length:lemma} 
we have $\dot{\lambda}(0)=
\ddot{\lambda} (0)=0$.
Also recall that $\dot{\eta}(0)=
\dot{\zeta} (0)=0$, so differentiating two times equality~(\ref{holcomp:eq}),
evaluating at $0$ and taking the trace we obtain
\begin{equation}\label{der:hol:eq}
\begin{array}{lllll}
2\ddot{\tr_{\overline{\gamma}}}(0) & = & \ddot{\eta}(0) & \cdot
& ((\lambda(0)+\lambda(0)^{-1})\cos(\eta (0)/2)\sin
(\zeta (0)/2)\\
& & & & + 2\sin(\eta(0)/2)\cos(\zeta (0)/2))\\
& + &
\ddot{\zeta}(0) & \cdot
& ((\lambda(0)+\lambda(0)^{-1})\sin(\eta (0)/2)\cos
(\zeta (0)/2)\\
& & & & + 2\cos(\eta(0)/2)\sin(\zeta (0)/2)).
\end{array}
\end{equation} 
Let us suppose 
$\ddot{\zeta}(0)=-\ddot{\eta} (0)\neq 0$. 
In this case,
from equation~(\ref{der:hol:eq}) we obtain
\[
\ddot{\tr_{\overline{\gamma}}}(0)=
(\ddot{\eta}(0)/2)(\lambda (0)+\lambda (0)^{-1}-2)
(\sin ((\zeta (0)/2)-(\eta (0)/2)))\neq 0.
\]

When
$\ddot{\zeta} (0)=0$ computations are more involved, and
in order to prove that $\ddot{\tr_{\overline{\gamma}}} (0)\neq 0$
one can show that
\begin{equation}\label{stima:eq}
(\lambda(0)+\lambda(0)^{-1})\cos(\eta (0)/2)\sin
(\zeta (0)/2) + 2\sin(\eta(0)/2)\cos(\zeta (0)/2)>0.
\end{equation}  
We skip this computation here, addressing the reader to~\cite{Fri:tesi}
for the details. The proof of Proposition~\ref{nonisol:prop}
is now concluded.

\section{Similar fillings}\label{similar:section}

Kojima proved in~\cite{Koj:proc} that every 
complete finite-volume hyperbolic manifold $N$ with non-empty
geodesic boundary admits a 
\emph{canonical decomposition} into geometric polyhedra. 
For later reference we record the following:

\begin{prop}\label{Koj:char:prop}
Any shortest
return path in $N$ is an edge of the Kojima decomposition of $N$.
Moreover, any compact regular partially truncated
tetrahedron isometrically immersed in $N$ whose 
internal edges are shortest return paths is a piece of the canonical 
decomposition of $N$.
\end{prop}

The following result is proved in~\cite{FriMaPe3}. 

\begin{teo} \label{main2:teo}
Let $M\in\calM_{g,k}$ with $\partial M=\Sigma_g 
\sqcup\big(\mathop{\sqcup}\limits_{i=1}^kT_i\big)$. Then the following holds:
\begin{enumerate}
\item \label{main:canonical:item} 
$M$ has a unique triangulation with $g+k$ tetrahedra, which gives
the canonical Kojima decomposition of $M$;
\item
The volume of the complete hyperbolic structure of
$M$ depends only on $g$ and $k$;
\item \label{main:genus:item}
The Heegaard genus of $\big(M,\Sigma_g,\mathop{\sqcup}\limits_{i=1}^kT_i\big)$ is $g+1$;
\item \label{main:homology:item}
$H_1(M;\matZ)=\matZ^{g+k}$;
\item \label{main:turaev:viro:item} 
The Turaev-Viro invariant $\TV_r(M)$ depends only on $r$, $g$ and $k$.
\end{enumerate}
\end{teo}

Manifolds in $\calM_{g,k}$ also share other geometric invariants:

\begin{teo}\label{aggiuntina:teo}
Let $M\in\calM_{g,k}$ be endowed with its complete hyperbolic structure. 
Then the following holds:
\begin{enumerate}
\item \label{cusp:volume:item}
The cusp volume of $M$ depends only on $g$ and $k$; 
\item \label{cusp:shape:item}
The Euclidean structures on the boundary tori of $M$ 
are all isometric to the regular hexagonal one;
\item \label{shortest:path:item}
The length of the shortest return path of $M$ 
depends only on $g$ and $k$.
\end{enumerate}
\end{teo}
\dimostraz
By Lemma~\ref{embedcusp:lemma},
a maximal regular horocusp neighbourhood for $M$ is obtained
by gluing the maximal horocusp neighbourhoods of the ideal vertices
of the non-compact tetrahedra of $\calT$, whence 
point~(\ref{cusp:volume:item}).
Point~(\ref{cusp:shape:item}) has already been established
and point~(\ref{shortest:path:item})
is a consequence of Propositon~\ref{Koj:char:prop} and
Theorem~\ref{main2:teo}-(\ref{main:canonical:item}).
\finedimo

\subsection{Spines and homology}
\begin{figure}
\begin{center}
\includegraphics[width=8cm]{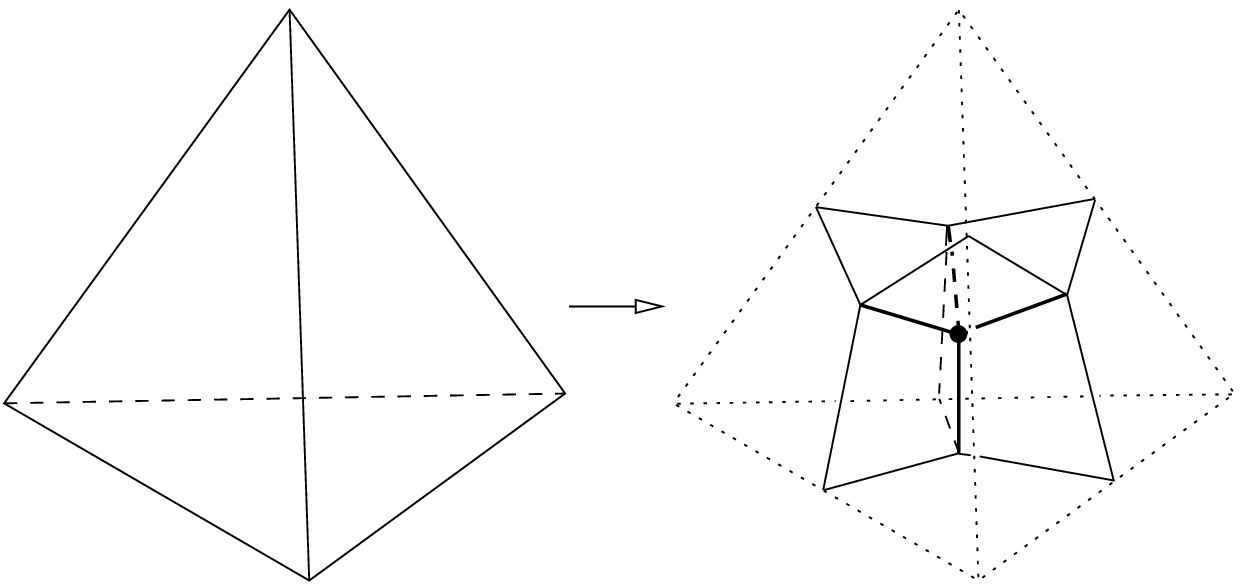}
\nota{From an ideal triangulation to a special spine.} \label{dualspine:fig}
\end{center}
\end{figure}

We now prove a refinement of 
Theorem~\ref{main2:teo}-(\ref{main:homology:item})
that will be useful later.
To this aim
we switch from the viewpoint of ideal
triangulations to the dual viewpoint of \emph{special spines},
suggested in Fig.~\ref{dualspine:fig}. Recall that a spine
of a manifold is a subpolyhedron onto which the manifold collapses.
A polyhedron is \emph{special} if it is locally homeomorphic to that of
Fig.~\ref{dualspine:fig}-right and its natural stratification consists
of $0$-, $1$-, and $2$-cells.

\begin{prop}\label{homology:prop}
Let $M\in\calM_{g,k}$. Then we have the
exact sequence 
\[
0\longrightarrow 
H_1\left(\mathop{\sqcup}\limits_{i=1}^k T_i;\matZ\right)
\mathop{\longrightarrow}\limits^{i_{\ast}}  
H_1(M;\matZ)\longrightarrow \matZ^{g-k}\longrightarrow 0,
\]
where $i_\ast$ is the map induced by
the inclusion $i:\mathop{\sqcup}\limits_{i=1}^k T_i
\to M$.
\end{prop}
\dimostraz
Let $P$ be the spine dual to the triangulation
of $M$ with $g+k$ tetrahedra. Note that $P$ has 
a cellularization into
vertices, edges, and faces corresponding to tetrahedra, faces, and edges of 
the triangulation. 
We denote in particular by $S(P)$ the 1-skeleton of $P$ (a 4-valent graph).
By Proposition~\ref{forma:tria:prop} the spine $P$ contains
$k$ (open) hexagonal faces $E_1,\ldots, E_k$ and one big face $G$ with $6g$
vertices (with multiplicity). 
For $i=1,\ldots,k$ the closure $\overline{E}_i$ of $E_i$ is a torus
which bounds a collar of the $i$-th toric component $T_i$ of $\partial M$,
and the rest of $P$ lies outside this collar.
Since $M$ collapses onto $P$,
we have $H_1(M;\matZ)\cong H_1(P;\matZ)$, and
we can use cellular homology to compute $H_1(P;\matZ)$.

Since $g>k$, a vertex $v$ of $P$ exists which
corresponds to a partially truncated tetrahedron
without ideal vertices. Notice that $G$ is the only face 
incident to $v$. Moreover, an easy analysis of the local
structure of $P$ near $v$ shows that the number of edges
emanating from $v$ on which $G$ passes three times with the same
orientation is at most two.
Let us also observe that any finite graph
has an even number of vertices with odd valence,
so the number of connected components of $S(P)\setminus
\{v\}$ is at most two. These facts easily imply that a maximal
tree $Y$ in $S(P)$ exists with the following properties:
$S(P)\setminus Y$ consists of $g+k+1$ edges
$e_1,\ldots,e_{g+k+1}$, where $e_{2i-1}$ and $e_{2i}$ 
represent generators of $H_1 (T_i;\matZ)$ for $i=1,\ldots,k$,
and $G$ passes three times on $e_{g+k+1}$ with different 
orientations.
Therefore we get a presentation for $H_1 (P;\matZ)$  
with generators $e_1,\ldots,e_{g+k+1}$ and 
one relator $w$ containing $\pm e_{g+k+1}$ once. 
This implies in turn that the classes of 
$e_1,\ldots,e_{g+k}$ give a free basis of
$\matZ^{g+k+1}/\langle w\rangle\cong\matZ^{g+k}\cong
H_1 (P;\matZ)$, whence the conclusion. 
\finedimo

\subsection{Boundary slopes and Dehn filling}\label{slopes:subsection}
Let $T_i$ be the $i$-th boundary torus of a manifold $M\in\calM_{g,k}$,
and recall that the unique complete finite-volume hyperbolic structure
on $M$ induces on $T_i$ a Euclidean structure defined up to similarity.
For the sake of simplicity, we endow $T_i$ with a fixed Euclidean structure
choosing the scale factor in
such a way that ${\rm Area}(T_i)=\sqrt{3}/2$. 
This easily implies that $T_i$ is isometric
to $\matC/\Gamma$, where $\matC$ is endowed with
the standard Euclidean metric, and $\Gamma$ is the 
discrete additive subgroup of $\matC$ with generators
$1,-1/2+i\sqrt{3}/2$. 
We denote by 
$\calM (T_i)$ the group of isotopy classes of isometries
of $T_i$.
Of course $\calM (T_i)$ acts
on the set of slopes on $T_i$.

Let
$D_6$ be the 
dihedral group with $12$ elements,
\emph{i.~e.}~the group of isometries of $\matC$
generated by the rotation $r:\matC\to\matC,
\ r(z)=e^{i\pi/3}\cdot z$ and the reflection
$s:\matC\to\matC,\ s(z)=\overline{z}$.
Any element of $D_6$ induces an isometry of $T_i$,
and any isometry of $T_i$ lifts up to isotopy to an element
of $D_6$. Thus $\calM (T_i)$ is canonically isomorphic to $D_6$.

Let $\mu_i,\lambda_i$ be the preferred basis of
$H_1 (T_i;\matZ)$ chosen in Subsection~\ref{compl:subsection}.
In what follows we will often
represent slopes as indivisible elements in $H_1(T_i;\matZ)$
without emphasizing the fact that each slope corresponds
in fact to two such elements.
Any slope $s$ on $T_i$ determines a well-defined isotopy class of geodesics
on $T_i$, and we denote by $L(s)$ the Euclidean length of such geodesics.
An elementary calculation gives the following:

\begin{lemma}\label{length:slope:lem}
Let $s=p\cdot \mu_i+ q\cdot \lambda_i$ be a slope on $T_i$. 
Then $L(s)=\sqrt{p^2+q^2-pq}$.
\end{lemma}

Let $\{\kappa_1<\kappa_2<\ldots<\kappa_n
<\ldots\}$ be the set of lengths of slopes on $T_i$.
The following result is easily deduced from Lemma~\ref{length:slope:lem}.

\begin{prop}\label{short:slopes:prop}
The following holds:
\begin{itemize}
\item
There are exactly three slopes of length $\kappa_1=1$.
They are represented by $\mu_i, \lambda_i$ and $\mu_i+\lambda_i$, and they
are $\calM (T_i)$-equivalent to each other.
\item
There are exactly three slopes of length $\kappa_2=\sqrt{3}$. 
They are represented by $\mu_i-\lambda_i, \mu_i+2\cdot \lambda_i$ 
and $2\cdot \mu_i+\lambda_i$, and they
are $\calM (T_i)$-equivalent to each other.
\item
If $s$ is a slope with $L(s)\geqslant \kappa_3=\sqrt{7}$, then there exist
exactly six slopes 
$\calM (T_i)$-equivalent to $s$.
\end{itemize}
\end{prop}

\begin{rem}\label{diff:length:rem}
Let $s$ and $s'$ be slopes on $T_i$.
Of course if $s'$ is $\calM (T_i)$-equivalent
to $s$ then $L(s')=L(s)$, but the converse
is not true. For example, the slopes $s=19\mu_i+11 \lambda_i$ and
$s'=16\mu_i-\lambda_i$ have the same length $L(s)=L(s')=
\sqrt{273}$, even if they are not $\calM (T_i)$-equivalent. 
\end{rem}

\subsection{Dehn fillings}
The following result, which is proved in~\cite{FriMaPe3}, completely classifies
the Dehn fillings of elements in $\calM_{g,k}$.

\begin{teo}\label{main:Dehn:teo}
Let $M\in\calM_{g,k}$ with 
$\partial M = \Sigma_g\sqcup\big(\mathop{\sqcup}\limits_{i=1}^k T_i\big)$,
let $h\leqslant k$,
let $s_i$ be a slope on $T_i$ for $i=1,\ldots, h$
and $N=M(s_1,\ldots,s_h)$. Then $N$ is hyperbolic if and only
if $L (s_i)\geqslant \kappa_3$ for all $i=1,\ldots,h$. Moreover,
when $N$ is hyperbolic the Heegaard genus of 
$(N,\Sigma_g,T_{h+1}\sqcup\ldots\sqcup T_{k})$
is $g+1$.
\end{teo}

Actually, it is proved in~\cite{FriMaPe3} that each manifold in $\calM_{g,k}$
is a link complement in the handlebody of genus $g$.
The following result is an easy consequence of
Proposition~\ref{homology:prop}.

\begin{prop}\label{hom:filling:prop}
If 
$N$ is as in the statement of Theorem~\ref{main:Dehn:teo}, then
$H_1 (N;\matZ)\cong \matZ^{g+k-h}$.
\end{prop}

\subsection{Symmetries of $\Omega_{g,k}$}
We now describe the symmetries of the deformation space
$\Omega_{g,k}$, and explain how these symmetries act on the space of
Dehn filling coefficients. 
In order to clarify our arguments it is convenient to 
denote the coordinates of $\matR^{12k+1}$ as in equation~(\ref{xxxx:eq}):
\[ 
\beta (x)=x_{12k+1},
\, \alpha_l^j (x)=x_{6(l-1)+j},\, \gamma_l^j (x)=x_{6(l-1)+3+j},
\ l=1,\ldots,2k,\, j=1,2,3.
\]

\subsection{Symmetries of $\Omega_{g,k}$}
Let ${\rm Aut}
(\Omega_{g,k})$ denote the set of diffeomorphisms of
$\Omega_{g,k}$ onto itself, fix $i\in\{1,\ldots,k\}$
and take an element $\sigma$ of the symmetric group
$\mathfrak{S}_3$. We can make $\sigma$ act 
on the apices of the dihedral angles of $\Delta_{2i-1}$ and of $\Delta_{2i}$,
thus obtaining 
an automorphism 
$\widehat{\sigma}_i\in{\rm Aut}(\Omega_{g,k})$ which
leaves the angles of all the other 
tetrahedra unchanged: 
\begin{equation}\label{def:sig:eq}
\begin{array}{lllllll}
\alpha_l^j (\widehat{\sigma}_i (x))&=&\alpha_l^{\sigma^{-1} (j)} (x),&
\gamma_l^j (\widehat{\sigma}_i (x))&=&\gamma_l^{\sigma^{-1} (j)} (x)&
{\rm if}\ l=2i-1,2i;\\
\alpha_l^j (\widehat{\sigma}_i (x))&=&\alpha_l^j (x),&
\gamma_l^j (\widehat{\sigma}_i (x))&=&\gamma_l^j (x)
& {\rm if}\ l\neq2i-1,2i;\\
\beta (\widehat{\sigma}_i (x))&=&\beta (x).&&&&
\end{array}
\end{equation}
The fact that $\widehat{\sigma}_i$ takes indeed $\Omega_{g,k}$ into
itself is a consequence of the invariance of consistency 
equations~(\ref{length:edge:eq}),
(\ref{sigma:eq}), (\ref{sum:eq})
under the permutation of apices described in~(\ref{def:sig:eq}).

Another symmetry $\zeta_i:\Omega_{g,k}\to\Omega_{g,k}$
exists which corresponds to interchanging the r\^oles of the
tetrahedra $\Delta_{2i-1}$ and $\Delta_{2i}$:
\begin{equation}\label{def:zeta:eq}
\begin{array}{lll}
\beta (\zeta_i (x))=\beta (x); &
\alpha_{2i-1}^j (\zeta_i (x))=\alpha_{2i}^{j} (x),&
\gamma_{2i-1}^j (\zeta_i (x))=\gamma_{2i}^{j} (x);\\
\alpha_l^j (\zeta_i (x))=\alpha_l^j (x),&
\gamma_l^j (\zeta_i (x))=\gamma_l^j (x)
& {\rm if}\ l\neq2i-1,2i.
\end{array}
\end{equation}
Also in this case the fact that $\zeta_i (\Omega_{g,k})=\Omega_{g,k}$
easily follows from a straight-forward analysis of the consistency
equations.
We can now define a map
\begin{equation}\label{phii:eq}
\varphi_i:\mathfrak{S}_3\times \matZ/_2\to{\rm Aut}(\Omega_{g,k}),
\
\varphi_i (\sigma,\epsilon)=\widehat{\sigma}_i\circ \zeta_i^\epsilon,
\quad \sigma\in \mathfrak{S}_3,\ \epsilon=0,1.
\end{equation}
Using that $\zeta_i$ commutes with $\widehat{\sigma}_i$ for all 
$\sigma\in \mathfrak{S}_3$
it is easily seen that $\varphi_i$ is an injective homomorphism
with image a certain subgroup ${\rm Sym}_i (\Omega_{g,k})$
of ${\rm Aut}(\Omega_{g,k})$.

If $\kappa$ is an element of the symmetric group 
$\mathfrak{S}_k$, then there exists a symmetry  
$\widetilde{\kappa}\in{\rm Aut}(\Omega_{g,k})$ which induces 
the corresponding permutation of the shape of the cusps:
\begin{equation}\label{def:kappa:eq}
\begin{array}{lll}
\alpha_{2i-1}^j (\widehat{\kappa} (x))=\alpha^j_{2\kappa^{-1} (i)-1} (x),&
\alpha_{2i}^j (\widehat{\kappa} (x))=\alpha^j_{2\kappa^{-1} (i)} (x),&\\
\gamma_{2i-1}^j (\widehat{\kappa} (x))=\gamma^j_{2\kappa^{-1} (i)-1} (x),&
\gamma_{2i}^j (\widehat{\kappa} (x))=\gamma^j_{2\kappa^{-1} (i)} (x),&
i=1,\ldots,k,\, j=1,2,3;\\
\beta(\widehat{\kappa}(x))=\beta (x). & &
\end{array}
\end{equation}
The map $\kappa\mapsto\widehat{\kappa}$ defines an injective 
homomorphism
$\nu:\mathfrak{S}_k\to {\rm Aut}(\Omega_{g,k})$.

Let us denote by
${\rm Sym} (\Omega_{g,k})$ the subgroup of 
${\rm Aut}(\Omega_{g,k})$ generated by $\nu(\mathfrak{S}_k)\cup \left(
\bigcup_{i=1}^k {\rm Sym}_i (\Omega_{g,k})\right)$.
Elements in ${\rm Sym}_i (\Omega_{g,k})$
commute with elements in ${\rm Sym}_j (\Omega_{g,k})$
whenever $i\neq j$, thus the group generated by
the ${\rm Sym}_i (\Omega_{g,k})$'s is actually isomorphic
to the product $\prod_{i=1}^k {\rm Sym}_i (\Omega_{g,k})$.
Moreover if 
$\pi:{\rm Sym}(\Omega_{g,k})\to \mathfrak{S}_k$ is the natural homomorphism
which maps each symmetry
to the corresponding permutation of cusps we have
${\rm Ker}\ \pi=\prod_{i=1}^k {\rm Sym}_i (\Omega_{g,k})$
and $\pi\circ\nu={\rm Id}:\mathfrak{S}_k\to \mathfrak{S}_k$. Thus
\begin{equation}\label{deco1:eq}
{\rm Sym} (\Omega_{g,k})=\nu (\mathfrak{S}_k)\ltimes  
\left(\prod_{i=1}^{k} 
{\rm Sym}_i (\Omega_{g,k})\right)
\cong \mathfrak{S}_k\ltimes\left(\prod_{i=1}^{k} 
\mathfrak{S}_3\times \matZ/_2\right).
\end{equation}
Our next task is to investigate how symmetries in ${\rm Sym}(\Omega_{g,k})$
act on the space of Dehn filling coefficients parameterizing  
a small neighbourhood of $x_0$
in $\Omega_{g,k}$.

\subsection{Action on Dehn filling coefficients}
Let us denote by
$\calM (T_1 \sqcup\ldots\sqcup T_k )$ the group of isotopy
classes of isometries of $T_1  \sqcup\ldots\sqcup T_k $.
For $\sigma\in\mathfrak{S}_k$ we define an
element $\eta (\sigma)\in \calM (T_1 \sqcup\ldots\sqcup T_k )$
which permutes the \emph{marked} tori $T_1,\ldots,T_k$
according to $\sigma$. Namely,
$\eta (\sigma)$ is the isotopy class of 
any element $\sigma'\in{\rm Isom}
(T_1\sqcup\ldots\sqcup T_k)$
with the following properties:
$\sigma' (T_j)=T_{\sigma (j)}$,
$\sigma'_\ast (\mu_j)=\mu_{\sigma (j)}$,
$\sigma'_\ast (\lambda_j)=\lambda_{\sigma (j)}$
for $j=1,\ldots,k$.
It is easily seen that the map 
$\eta:\mathfrak{S}_k\to
\calM (T_1\sqcup\ldots\sqcup T_k)$ is a well-defined
injective homomorphism.
Let now $\pi':\calM (T_1\sqcup\ldots\sqcup T_k)\to
\mathfrak{S}_k$ be the natural projection which
associates to any element in $\calM (T_1\sqcup\ldots\sqcup T_k)$
the induced permutation of the $T_j$'s.
Of course we have $\pi'\circ \eta={\rm Id}:
\mathfrak{S}_k\to\mathfrak{S}_k$, and the kernel
of $\pi'$ is canonically isomorphic to
$\calM (T_1)\times\ldots\times \calM (T_k)$. Therefore
we have
\begin{equation}\label{deco2:eq}
\calM(T_1\sqcup\ldots\sqcup T_k)=
\eta(\mathfrak{S}_k)\ltimes \left(\prod_{i=1}^{k} 
\calM (T_i)\right)
\cong \mathfrak{S}_k\ltimes\left(\prod_{i=1}^{k} 
D_6\right),
\end{equation}
where $D_6\cong \mathfrak{S}_3\times \matZ/_2$ is the dihedral
group with $12$ elements.

By Theorems~\ref{smoothness:teo} and~\ref{essential:bis:teo},
from now on we can fix
a small neighbourhood $V$ of $x_0$ in $\Omega_{g,k}$
such that
for all $x\in V$ the Dehn filling coeffcient
$(p_j(x),q_j(x))\in S^2=\matR^2\cup\{\infty\}$ 
is well-defined, and the map
\begin{equation}\label{definingd:eq}
d=(d_1,\ldots,d_k):V\to \prod_{i=1}^k S^2,\quad
d_j (x)=(p_j(x),q_j(x))\in S^2
\end{equation} 
is a diffeomorphism onto an open 
neighbourhood
of $\{\infty\}\times\ldots\times\{\infty\}$
in $S^2\times\ldots\times S^2$.
It is easily seen that we can also assume
$\psi(V)=V$ for all $\psi\in{\rm Sym}(\Omega_{g,k})$.

We now observe that any element $h\in\calM(T_1\sqcup\ldots \sqcup T_k)$ 
induces an automorphism
of $H_1 (T_1;\matR)\oplus\ldots\oplus H_1 (T_k;\matR)$.
The basis $\mu_i,\lambda_i$ defines a canonical
isomorphism $H_1 (T_i;\matR)\cong\matR^2$, so $h$ induces
an automorphism $h_\ast$
of $\prod_{i=1}^{k} S^2$ that preserves $\{\infty\}\times\ldots\times
\{\infty\}$.

\begin{prop}\label{iso:to:iso:prop}
For any $\psi\in{\rm Sym}(\Omega_{g,k})$ there exists a unique
$h(\psi)\in\calM(T_1\sqcup\ldots\sqcup T_k)$ such that
$d(\psi(x))=h(\psi)_\ast (d(x))$ for all $x\in V$. 
Moreover the map
\[
{\rm Sym}(\Omega_{g,k})\to\calM (T_1\sqcup\ldots\sqcup T_k),
\quad \psi\mapsto h(\psi)
\]
is a group isomorphism which preserves decompositions~(\ref{deco1:eq}),
(\ref{deco2:eq}). 
\end{prop}
\dimostraz
We need to describe as explicitly as possible
how the maps $d_j:V\to S^2$, $i=1,\ldots,k$ change under
precompositions with elements in ${\rm Sym} (\Omega_{g,k})$.
Let us
consider the action
of $\mathfrak{S}_3\times\matZ/_2$ on $\Omega_{g,k}$ via the representation
$\varphi_j:\mathfrak{S}_3\times\matZ/_2\to{\rm Sym}_j (\Omega_{g,k})$
defined in equation~(\ref{phii:eq}).
By definition the element $(\sigma,0)\in\mathfrak{S}_3\times\matZ/_2$
acts on the $\gamma_{2j-1}^l$'s and the $\gamma_{2j}^l$'s just
by applying $\sigma^{-1}$ to the apices, while the action of 
$({\rm Id},1)\in\mathfrak{S}_3\times\matZ/_2$ 
interchanges the indices $2j-1, 2j$.
Let $r=((132),1)$, $s=((12),0)$
be fixed elements of $\mathfrak{S}_3\times\matZ/_2$.
Together with equations~(\ref{uv:eq}),
the description of the action of $\varphi_j (r),
\varphi_j (s)$ given above implies (after some computations)
that
\begin{eqnarray}
u_j (\varphi_j (r)(x))
= -v_j (x),&
v_j (\varphi_j (r)(x))
= u_j (x)+v_j (x)\label{newuv1:eq},\\
u_j(\varphi_j (s)(x)) 
= -\overline{u_j (x)},&
v_j (\varphi_j (s)(x))
=\overline{u_j (x)}+\overline{v_j(x)}.\label{newuv2:eq} 
\end{eqnarray}
We can now compute the action of $h(\varphi_j (r))$
and $h(\varphi_j (s))$ on Dehn filling coefficients.
Using equations~(\ref{newuv1:eq}),~(\ref{newuv2:eq}) and the very definition of
Dehn filling coefficients we get
\[
\begin{array}{lll}
2\pi i &=& 
p_j (\varphi_j (r)(x))u_j (\varphi_j (r)(x))
+q_j (\varphi_j (r)(x)) v_j (\varphi_j (r)(x))\\ &=&
p_j (\varphi_j (r)(x)) (-v_j (x))
+q_j (\varphi_j (r)(x)) (u_j (x)+v_j(x))\\ &=&
q_j (\varphi_j (r) (x)) u_j (x)+
(q_j (\varphi_j (r) (x))-p_j (\varphi_j (r) (x))) v_j (x),
\end{array}
\]
whence $p_j (x)=q_j (\varphi_j (r) (x))$, 
$q_j (x)=q_j (\varphi_j (r) (x))-p_j (\varphi_j (r) (x))$ and
\begin{equation}\label{primautile:eq}
p_j (\varphi_j (r)(x))=p_j (x)- q_j (x),\quad
q_j (\varphi_j (r) (x))=p_j (x).
\end{equation}
A similar computation also gives
\begin{equation}
p_j (\varphi_j (s)(x))=p_j (x)- q_j (x),\quad
q_j (\varphi_j (s) (x))=-q_j (x). 
\end{equation} 
This easily implies that $h(\varphi_j (r))$ and $h(\varphi_j (r))$
act on Dehn filling coefficients at the $j$-th end of $M$ respectively
as a rotation of angle $\pi/3$ and a reflection with respect
to the line $\matR\cdot\mu_j$.
This gives in turn that $h$ restricts to an isomorphism
\[
{\rm Sym}(\Omega_{g,k})\supset {\rm Sym}_j (\Omega_{g,k})
\cong
\calM (T_j)\subset \calM(T_1\sqcup\ldots\sqcup T_k).
\]
Moreover, with notation as in formulae~(\ref{deco1:eq}),
(\ref{deco2:eq}), any permutation of cusps in 
$\nu (\mathfrak{S}_k)\subset {\rm Sym}(\Omega_{g,k})$ 
is taken by $h$ into the corresponding permutation in
$\eta(\mathfrak{S}_k)\subset \calM (T_1\sqcup\ldots\sqcup T_k)$, and this
concludes the proof.
\finedimo

\subsection{Return paths}
Recall that for $x\in\Omega_{g,k}$ we denote by $M(x)$
the hyperbolic structure induced on $M$ by $x$, and by
$\widehat{M} (x)$ the metric completion of 
$M(x)$. Moreover, if $x\in I\Omega_{g,k}$ then $\widehat{M} (x)$
is a complete finite-volume hyperbolic manifold with
geodesic boundary. In this case the unique compact edge in the
geometric triangulation of $M(x)$ defines a return path $l_x$
in $\widehat{M} (x)$. For $y\in I\Omega_{g,k}$
we denote by $L^y$ the length
with respect to the hyperbolic metric on $\widehat{M} (y)$.
Of course we have
$\lim_{y\in I\Omega_{g,k},\, y\to x_0}
L^y (l_y)=L^{x_0} (l_{x_0})$. Moreover,
a positive number $\delta$ exists such that
any return path in $M$ different from
$l_{x_0}$ has length at least $L^{x_0} (l_{x_0})+2\delta$.

\begin{lemma}\label{short:return:lemma}
There exists a neighbourhood $V'\subset V$ of
$x_0$ in $\Omega_{g,k}$ such that if
$x\in I\Omega_{g,k}\cap V'$, then $l_x$ is the only
return path in $\widehat{M} (x)$ having
length less than $L^{x_0} (l_{x_0})+\delta$.
In particular, if
$x\in I\Omega_{g,k}\cap V'$ then $l_x$ is the unique
shortest return path in $\widehat{M} (x)$.
\end{lemma}
\dimostraz
We suppose by contradiction that there exists a sequence
$\{y_n\}_{n\in\matN}\subset I\Omega_{g,k}$ converging to $x_0$
such that $\widehat{M} (y_n)$ contains a return path
$l_{n}\neq l_{y_n}$ with $L^{y_n} (l_n)< L^{x_0} (l_{x_0})+\delta$
for every
$n\in\matN$. 
Since the distance between the added geodesics
$\widehat{M}(y)\setminus M(y)$ and the geodesic boundary of
$\widehat{M}(y)$ approaches $\infty$ as $y$ tends to $x_0$,
we can suppose that the compact set $K_y\subset \widehat{M} (y)$ 
of points whose distance from $\partial \widehat{M} (y)$ is 
less than or equal to $2 L^{x_0}  (l_{x_0})$ is entirely 
contained in $M (y)$. 
Moreover, up to passing to a subsequence
we can suppose
that there exists an $\epsilon_n$-biLipschitz homeomorphism
$f_n: K_{y_n}\to K_{x_0}$ taking $\partial \widehat{M} (y_n)$
onto $\partial \widehat{M} (x_0)=\partial M (x_0)$, 
where $\epsilon_n$ tends to $1$
as $n$ tends to $\infty$.
Thus $L^{x_0} (f_{N} (l_{N}))< L^{x_0} (l_{x_0})+2\delta$ for some
$N \gg 0$. Since $l_{N}$ is not 
boundary-parallel  in $\widehat{M} (y_n)$, the path $f_{N} (l_{N})$
is not boundary-parallel in $M(x_0)$. 
Since return paths minimize length in their relative homotopy class,
this implies that $f_{N} (l_{N})$
is homotopic to $l_{x_0}=f_{N} (l_{y_{N}})$
relatively to the boundary in $M(x_0)$, so $l_{N}$
is homotopic to $l_{y_{N}}$ relatively to the boundary in
$\widehat{M} (y_n)$. Since $l_N$ and $l_{y_n}$
are both return paths, this implies in turn
$l_{N}=l_{y_{N}}$, a contradiction.
\finedimo
    
Let $V'$ be a neighbourhood of $x_0$ in $\Omega_{g,k}$
as in the statement of Lemma~\ref{short:return:lemma}, and
for $x\in I\Omega_{g,k}\cap V$ let $U(x)$ be the universal covering
of the hyperbolic manifold $\widehat{M} (x)$. 
We recall that $U(x)$ is isometric to a convex polyhedron of $\matH^3$
bounded by a countable number of hyperbolic planes.
Lemma~\ref{short:return:lemma} readily implies the following:

\begin{cor}\label{short:return2:cor}
The minimal distance between distinct connected components
of $\partial U(x)$ is equal to $L^x (l_x)$ for all
$x\in I\Omega_{g,k}\cap V'$. Moreover, if $S_1,S_2$ are distinct connected
components of $\partial U(x)$, then the distance between $S_1$ and $S_2$
equals $L^x (l_x)$ if and only if the shortest path joining $S_1$
with $S_2$ projects onto $l_x\subset \widehat{M} (x)$.
\end{cor}

The following proposition relates the intrinsic geometry of $U(x)$
to properties of our geometric triangulation
of $M(x)$, when $x\in I\Omega_{g,k}$.

\begin{prop}\label{univ:intrinsic:prop}
There exists a neighbourhood $V''\subset V'$ of
$x_0$ in $\Omega_{g,k}$ with the following property. 
Let $x\in I\Omega_{g,k}\cap V''$ and $S_1,\ldots,S_4$ be pairwise
distinct connected
components of $\partial U(x)$. Then the distance between $S_i$ and 
$S_j$ equals $L^x (l_x)$ for all $i\neq j$ if and only if
there exists a lift of a compact tetrahedron in the geometric
triangulation parameterized by $x$ whose truncation triangles
lie on $S_1,\ldots, S_4$.
\end{prop}
\dimostraz
We concentrate on the ``only if'' part of the statement,
the ``if'' part being obvious.
Let $\delta$ be as in the statement of Lemma~\ref{short:return:lemma}.
Then there exist $\varepsilon>0$ and a small neighbourhood $V''\subset
V'$ of $x_0$ in $\Omega_{g,k}$ such that for every $y\in I\Omega_{g,k}
\cap V''$ we have
\[
(1+\varepsilon) L^y (l_y)<L^{x_0}(l_{x_0})+\delta.
\]
Let $K_y$ be the set of points  of $\widehat{M}(y)$ whose distance
from $\partial \widehat{M} (y)$ is at most twice the diameter
of the regular truncated tetrahedron with edge-length equal
to $L^y (l_y)$. Up to resizing $V''$ we can suppose that
for all $y\in V''$ the set
$K_y$ is contained in $M(y)\subset \widehat{M} (y)$, and
there exists a 
$(1+\varepsilon)$-biLipschitz homeomorphism
$p_y:K_{y}\to K_{x_0}$. 

Let now $x\in I\Omega_{g,k}\cap V''$ and $S^x_1,\ldots,S^x_4$ be
pairwise
distinct connected
components of $\partial U(x)$ such that the distance between $S^x_i$ and 
$S^x_j$ equals $L^x (l_x)$ for all $i\neq j$. Let $\Delta$ be a topological
partially truncated tetrahedron with truncation triangles $B_1,B_2,B_3,B_4$
and internal edges $e_{ij}$ joining $B_i$ with $B_j$. 
We consider a geometric realization
$\widetilde{r}_x:\Delta\hookrightarrow U(x)\subset\matH^3$ with 
$\widetilde{r}_x (B_i)\subset S_i^x$, $i=1,2,3,4$, and we observe that $\widetilde{r}_x
(e_{ij})$
is the
shortest geodesic arc joining $S_i^x$ with $S_j^x$. 
Let $r_x: \Delta\to \widehat{M} (x)$ be the composition of
$\widetilde{r}_x$ with the projection $U(x)\to\widehat{M} (x)$.
Since $r_x(\Delta)\subset
K_x$, we can consider the map 
$r_{x_0}= p_x\circ r_x:\Delta\to \widehat{M}(x_0)=M (x_0)$, which lifts
in turn to $\widetilde{r}_{x_0}:\Delta\to U(x_0)$.
For  $i=1,2,3,4$
let $S^{x_0}_i$ be the component of $\partial U (x_0)$ containing
$\widetilde{r}_{x_0} (B_i)$, 
so that $\widetilde{r}_{x_0} (e_{ij})$ is a (not necessarily geodesic)
arc joining $S^{x_0}_i$ with $S^{x_0}_j$.
Since $r_x (e_{ij})$ is not null-homotopic relatively to the boundary in 
$\widehat{M} (x)$, $r_{x_0} (e_{ij})$ is not null-homotopic relatively
to the boundary in $U (x_0)$, so
$S^{x_0}_i\neq S^{x_0}_j$ for  $i\neq j$. Moreover 
we have $L^{x_0} (r_{x_0} (e_{ij}))\leqslant (1+\varepsilon) L^x (l_x)$,
whence
\[
d_{x_0} (S^{x_0}_i,S^{x_0}_j)\leqslant(1+\varepsilon) L^{x} (l_x) <
L^{x_0} (l_{x_0})+\delta,\quad i,j\in\{1,2,3,4\},\, i\neq j,
\]
where we denote by $d_{x_0}$ the hyperbolic metric on $U(x_0)$.
Thus
by Lemma~\ref{short:return:lemma} and Proposition~\ref{Koj:char:prop}
the hyperbolic planes $S^{x_0}_1,S^{x_0}_2,
S^{x_0}_3,S^{x_0}_4$ bound a compact 
geodesic regular partially truncated tetrahedron which 
projects onto a piece of the Kojima
decomposition of $M(x_0)$. This easily implies that 
$\widetilde{r}_x (\Delta)$ projects onto a compact partially truncated tetrahedron
in the geometric triangulation of $M(x)$.
\finedimo

\subsection{Similar fillings}\label{similar2:section}
A \emph{set of slopes}
for a complete finite-volume hyperbolic $3$-manifold $N$
is a set $\calS=\{s_{i_1},\ldots,s_{i_h}\}$ of either
$0$ or $1$ slope per boundary torus.
If $\calS=\{s_{i_1},\ldots,s_{i_h}\}$ is a set of slopes for $N$
we denote by $N(\calS)$ the manifold obtained by filling
$N$ along $s_{i_1},\ldots,s_{i_h}$.

Let $M$ and $M'$ be elements in 
$\calM_{g,k}$ (we do not exclude the case $M=M'$)
with boundary tori 
$T_1,\ldots,T_k$ and $T'_1,\ldots,T'_k$. 
We endow each of these tori
with the Euclidean metric defined on them by the hyperbolic
structures on $M, M'$ together with the requirement
that ${\rm Area} (T_i)={\rm Area} (T'_i)=\sqrt{3}/2$ for $i=1,\ldots,k$. 
We say that a set of slopes $\calS$ for $M$ is \emph{equivalent}
to the set of slopes $\calS'$ for $M'$
if there exists an orientation-preserving isometry
$\psi:T_1\sqcup\ldots\sqcup T_k\to T'_1\sqcup\ldots\sqcup T'_k$
taking $\calS$ onto $\calS'$. Of course if $\calS$ is equivalent
to $\calS'$ then the lengths
of the slopes in $\calS$ are equal to the lengths of the slopes in $\calS'$.
We recall however that the converse is not true (see Remark~\ref{diff:length:rem}).

\begin{teo}\label{similar:first:teo}
Let $M,M'$ be elements of $\calM_{g,k}$ and 
$\calS$ (resp.~$\calS'$) be a set of slopes for $M$ (resp.~$M'$).
Then there exists a positive constant $C$
such that the following holds:
if all the slopes of $\calS$ are longer than $C$ and 
$\calS$ is equivalent to $\calS'$, then
$M(\calS)$ is geometrically similar to $M' (\calS')$.
\end{teo}
\dimostraz
Let $V'$ be a neighbourhood of $x_0$
in $\Omega_{g,k}$ as in the statement of Lemma~\ref{short:return:lemma}
and $d:V'\to S^2\times\ldots\times S^2$ be the map defined
in equation~(\ref{definingd:eq}).
We can choose a positive constant $C$ depending only on $g$ and 
$k$ such that the following holds: if $\calS=\{s_{i_1},\ldots,
s_{i_h}\}$ is a set of slopes for $M$ with $L(s_{i_l})>C$
for $l=1,\ldots,h$, then any $k$-uple
of Dehn filling coefficients corresponding to $\calS$
lie in $d(V')$ (due to the choice of the signs, there exist exactly
$2^h$ such $k$-uples).

Let now $\calS$ be a set of slopes for $M$ whose elements are
longer than $C$ and let $\calS'$ be a set of slopes for $M'$
which is equivalent to $\calS$. Choose also points $x,x'\in V'\subset
\Omega_{g,k}$
such that $d(x)$ (resp.~$d(x')$) gives a $k$-uple of Dehn filling
coefficients corresponding to $\calS$ (resp.~$\calS'$).
By Proposition~\ref{iso:to:iso:prop}
it follows that a symmetry $\varphi\in {\rm Sym}(\Omega_{g,k})$ exists
with $\varphi (x)=x'$. Let $M(x)$ (resp.~$M'(x')$) be the hyperbolic
structure defined by $x$ on $M$ (resp.~by $x'$ on $M'$). 
Recall that $M(\calS)$ (resp.~$M'(\calS')$) 
is isometric to the metric completion of
$M(x)$ (resp.~of $M' (x')$), and that $M(\calS)\setminus M(x)$
(resp.~$M' (\calS')\setminus M' (x')$) is the union
of $h$ disjoint geodesics in $\calM(\calS)$ (resp.~in $\calM(\calS')$).
Notice
that since $x$ is ${\rm Sym} (\Omega_{g,k})$-equivalent to $x'$,
the geometric partially truncated tetrahedra in the decomposition of $M(x)$ 
are isometric to the geometric tetrahedra in the decomposition of $M'(x')$.
Together with Lemma~\ref{short:return:lemma}, this readily
implies that the shortest return paths of
$M(\calS)$ and $M(\calS')$ have the same length.
Moreover, $M(x)$ and $M'(x')$ have the same volume, whence
 ${\rm volume}(M(\calS))={\rm volume}(M' (\calS'))$.  

By Theorem~\ref{isolation:teo} the bases of
the cusps of $M(\calS)$ and $M' (\calS')$ are all isometric to
regular hexagonal tori, so $M(\calS)$
and $M' (\calS')$ share the same cusp shape.

Consider now the shape of the geometric tetrahedra in the 
triangulations $\calT=\{\Delta_1,\ldots,\Delta_{g+k}\}$, $\calT'=
\{\Delta'_1,\ldots,\Delta'_{g+k}\}$ of $M,M'$ respectively.
Without loss of generality we can order the tetrahedra of these triangulations
in such a way that $\Delta_l,\Delta'_l$ are asymptotic
to the cusps of $M(\calS), M' (\calS')$ for $l=2h+1,\ldots,2k$ (this
is equivalent to requiring that the slopes in $\calS$ and $\calS'$
lie on $T_1,\ldots,T_h$ and $T'_1,\ldots,T'_h$).
Then by Proposition~\ref{iso2:prop} 
a real number $\vartheta(x)=\vartheta (x')\in (0,\pi/3)$
exists
such that $\Delta_l$ and $\Delta'_l$ are isometric to
$\Delta^{\vartheta (x)}$ for $l=2h+1,\ldots, 2k$.
For $l=2h+1,\ldots,2k$ let now $v_l,v'_l$ be the ideal vertices
of $\Delta_l,\Delta'_l$ respectively. Due to Lemma~\ref{embedcusp:lemma}
and the symmetric
shape of $\Delta^{\vartheta (x)}$, up to increasing $C$
we can suppose that a unique
horocusp neighbourhood
$H_l$ 
of $v_l$ in $\Delta_l$ exists
which is tangent to the truncation triangles of
$\Delta_l$ and is entirely contained in
$\Delta_l$. Moreover $H_{2i-1}$ and $H_{2i}$ 
glue up in $M(\calS)$ giving a horocusp neighbourhood $O_i$
of the $i$-th cusp for $i=h+1,\ldots,k$.  
Also notice that the total horocusp neighbourhood
$O_{h+1}\sqcup\ldots\sqcup O_{k}$ is regular (since the $H_l$'s
are isometric to each other) and maximal (since each $O_i$ 
is tangent to the boundary of $M(\calS)$).
The very same construction also leads to a horocusp
neighbourhood $O'_{h+1}\sqcup\ldots\sqcup O'_{k}$ for $M'(\calS')$.
Since $\Delta_l$ is isometric to $\Delta'_l$ for
$l=2h+1,\ldots,2k$, we have ${\rm vol} (O_{h+1}\sqcup\ldots\sqcup O_{k})=
{\rm vol} (O'_{h+1}\sqcup\ldots\sqcup O'_{k})$, so
$M(\calS)$ and $M'(\calS')$ share the same cusp volume.

Recall now that Thurston's hyperbolic Dehn filling Theorem
ensures that if for all $l=1,\ldots,h$ we have
$L(s_{l})>C'>0$ for some sufficiently 
large $C'$, then
the shortest geodesics of $M(\calS)$ and $M'(\calS')$
are exactly the geodesics added to $M(x)$ and $M' (x')$.
Thus under the hypothesis that $L(s_{l})>C'$ for all
$l=1,\ldots,h$, in order to prove that the shortest geodesics
of $M(\calS)$ and $M'(\calS')$ have the same complex length 
we only have to compute the complex length of these added geodesics.
The desired result is then easily obtained from 
Proposition~\ref{added:geo:prop}
and equations~(\ref{newuv1:eq}),~(\ref{primautile:eq}).

The fact that $H_1 (M(\calS);\matZ)$ is isomorphic to 
$H_1 (M' (\calS');\matZ)$
is an immediate consequence of Proposition~\ref{hom:filling:prop}.
By Theorem~\ref{main:Dehn:teo},
if $\Sigma$ (resp.~$\Sigma'$) is the geodesic boundary
of $M(\calS)$ (resp.~of $M' (\calS')$), then both the Heegaard
genus of $(M(\calS),\Sigma,\partial M (\calS)\setminus \Sigma)$
and  the Heegaard
genus of $(M'(\calS'),\Sigma',\partial M' (\calS')\setminus \Sigma')$
are equal to $g+1$.

In order to prove our statement about Turaev-Viro invariants
we need to construct special spines for $M(\calS)$ and
$M'(\calS')$.
Let $P\subset M$ be the special spine of $M$ 
dual to the canonical decomposition $\calT$ of $M$, and
recall that for $j=1,\ldots,k$ a hexagonal face $E_j$
of $P$ exists which is parallel to the $j$-th boundary torus
of ${M}$. Let $E_{1},\ldots, E_{h}$ be the faces corresponding
to the filled tori in $M(\calS)$ and for $l=1,\ldots,h$
let $m_{l}$ be a loop on $E_{l}$ which represents the slope
$s_{l}\in\calS$ and is in general position with respect to the
singular locus $S(P)$ of $P$.
The complement 
of $P\subset M\subset M(\calS)$ inside
$M(\calS)$ consists of the disjoint union of 
an open collar of $\partial M(\calS)$ 
and $h$ open solid tori. 
Take meridinal discs $D_{1},\ldots,D_{h}$ of these solid tori
with $\partial D_{l}=m_{l}$ for $l=1,\ldots,h$.
The complement of $P\cup D_{1}\cup\ldots\cup D_{h}$ 
is as above, with $h$ open balls instead of
the $h$ open solid tori.
Fix now $l\in\{1,\ldots,h\}$.
The loop $\partial D_{l}\subset E_{l}$ cuts $E_{l}$ 
into several discal open faces of
$P\cup D_{1}\cup\ldots\cup D_{h}$ (these faces are
indeed homeomorphic to open discs because the loop $m_{l}$ is
sufficiently complicated with respect to the graph $S(P)\cap
\overline{E}_{l}$). 
Each such face
separates $\Sigma_g\subset\partial M(\calS)$ from the open ball
corresponding to the $l$-th added solid torus, 
so, if we remove from $P\cup D_{1}\cup
\ldots\cup D_{h}$ 
a face for each added solid torus we end up with
a special spine $P(\calS)$ of $M(\calS)$.
The very same procedure also provides a special spine
$P' (\calS')$ for $M'(\calS')$. Let $\calT(\calS),
\calT'(\calS')$ be the triangulations of $M(\calS),
M' (\calS')$ dual to $P(\calS), P'(\calS')$ respectively.
It is not difficult to show that since $\calS$ is equivalent
to $\calS'$, the loops representing the slopes in $\calS$ and
in $\calS'$ can be chosen so that the incidence numbers between
edges and tetrahedra are the same for $\calT(\calS)$
and for $\calT' (\calS')$. As pointed out 
in~\cite{MatNow}, this implies
that $M(\calS)$ and $M' (\calS')$ share the
same Turaev-Viro invariants. 

Let $\psi:T_1\sqcup\ldots \sqcup T_k\to
T'_1\sqcup\ldots \sqcup T'_k$ be the orientation-preserving
isometry taking $\calS$ onto $\calS'$, and $\calS_\ast$
be a set of sufficiently complicated slopes for $M(\calS)$.
Of course we can regard $\calS_\ast$ as a set of slopes
for $M$ too. Let $\calS'_\ast$ be the set of slopes for $M'$
obtained by applying $\psi$ to $\calS_\ast$. Of course
$\calS'_\ast$ is a set of slopes for $M' (\calS')$, and
the manifolds $M(\calS) (\calS_\ast)\cong M(\calS\cup\calS_\ast)$
and $ M'(\calS') (\calS'_\ast)\cong M'(\calS'\cup\calS'_\ast)$
share the same volume, homology, cusp volume, cusp shape, length
of the shortest return path, complex length of the shortest geodesic,
Heegaard genus 
and Turaev-Viro invariants.
\finedimo

\begin{rem}
Let $M,M'$ be elements of $\calM_{g,k}$ and suppose that
$s,s'$ are sufficiently long slopes on the tori $T\subset\partial {M},
T'\subset\partial {M}'$.
If there exists an \emph{orientation-reversing}
isometry of $T$ onto $T'$ taking $s$ into $s'$, then the complex length
of the added geodesic in $M(s)$ is equal to the conjugate of the
complex length of the added geodesic in $M'(s')$.
\end{rem}

\subsection{Non-homeomorphic fillings}\label{nonhomeo:subsection}
\begin{figure}
\begin{center}
\input{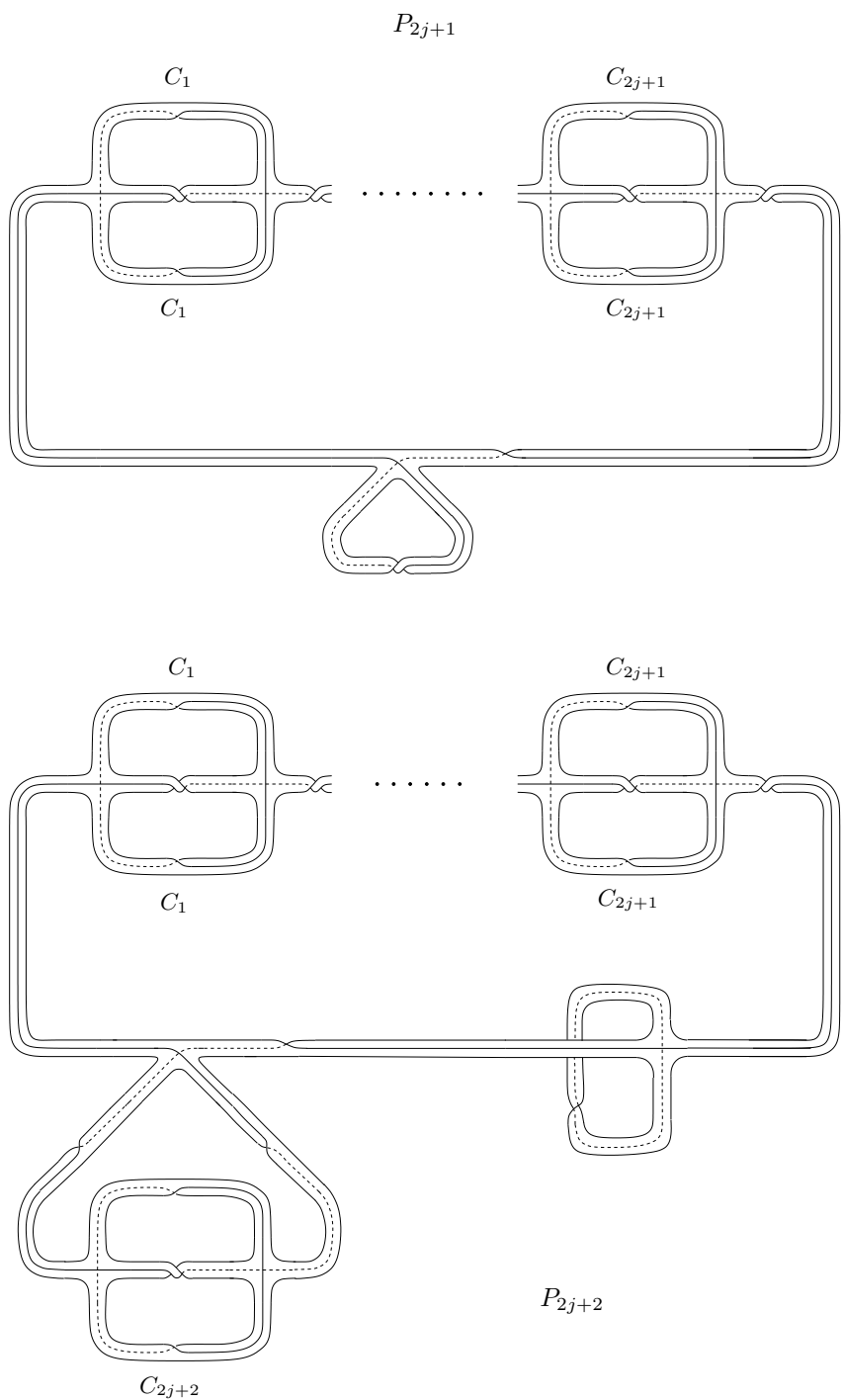}
\caption{{The regular neighbourhood 
of the $1$-skeleton $S(P_k)$ of $P_k$.
Each pair of vertices joined by three edges in $S(P_k)$
gives rise to a toric cusp in $X_k$.}}\label{casopar:fig}
\end{center}
\end{figure}

This paragraph is devoted to the proof of Theorem~\ref{nonhomeo:fil:teo}.
Let $P_k$ be the special polyhedron 
whose $1$-skeleton has the regular neighbourhood 
described in Fig.~\ref{casopar:fig}. It is easily seen that $P_k$ is the spine
of a manifold $X_k$. Computing the boundary of $X_k$ as explained 
in~\cite{BenPet:manuscripta} one can easily prove that $X_k\in
\calM_{k+1,k}$ if $k$ is odd and $X_k\in\calM_{k+2,k}$ if $k$
is even. 

\begin{prop}\label{iso:xk:prop}
For all $k\geqslant 1$, the manifold
$X_k$ admits no non-trivial isometries.
\end{prop}
\dimostraz
Let $\calT_k$ be the triangulation of $X_k$ dual to $P_k$.
Since $\calT_k$ is the Kojima decomposition of $X_k$,
the group of isometries of $X_k$ is canonically isomorphic
to the group
${\rm Aut} (\calT_k)$ of the combinatorial
automorphisms of $\calT_k$.
Now a straightforward analysis of the combinatorics of $\calT_k$
shows that ${\rm Aut} (\calT_k)$ is trivial, whence the conclusion.
\finedimo

Together with Proposition~\ref{iso:xk:prop}, the following result
implies Theorem~\ref{nonhomeo:fil:teo}.

\begin{prop}
Let $X\in\calM_{g,k}$ with boundary tori $T_1,\ldots,T_k$ and
suppose that $X$ admits no non-trivial isometry.
For each $i=1,\ldots,k$ we can choose a finite set $\calS_i$ of slopes on
$T_i$ with the following property. Let
$\calS$ be a set of slopes for $X$ whose elements do not belong to
$\calS_i$, $i=1,\ldots,k$ and let $h=\# \calS\leqslant k$.
Then the number of sets of slopes equivalent to
$\calS$ is greater than or equal to $(k!\cdot 3^h)/(h!\cdot (k-h)!)$. 
Moreover, if
$\calS'$ is a set of slopes equivalent to $\calS$ and $X (\calS)$
is homeomorphic to $X (\calS')$, then
$\calS=\calS'$.
\end{prop}
\dimostraz
Thurston's hyperbolic Dehn filling Theorem and Theorem~\ref{similar:first:teo}
imply that
we can choose the finite set $\calS_i$
in such a way that if $\calS$ is as in the statement
and $\calS'$ is a set of slopes equivalent to $\calS$, then
the following conditions hold:
no slope in $\calS'$ is contained in some $\calS_i$;
$X (\calS)$, $X (\calS')$
are geometrically similar hyperbolic $3$-manifolds; 
the cores of the added solid tori 
give the $h$ shortest geodesics both of $X (\calS)$
and of $X (\calS')$.

An elementary combinatorial argument shows that
the number of sets of slopes equivalent to $\calS$
is at least $(k!\cdot 3^h)/(h!\cdot (k-h)!)$.

Suppose now that $\calS'$ is equivalent to $\calS$ and let $\psi:X (\calS)\to
X (\calS')$ be a homeomorphism. By Mostow-Prasad's rigidity Theorem, 
$\psi$ is homotopic to an isometry $\psi'$, which  
must take the added geodesics of $X (\calS)$
to the added geodesics of $X (\calS')$. This gives in turn a homeomorphism
of $X$ onto itself taking $\calS$ onto $\calS'$.  
By rigidity again, up to homotopy such a homeomorphism
restricts to an isometry of $X$,  whence
$\calS=\calS'$ since $X$ admits no non-trivial isometry.
\finedimo

\section{Commensurability of similar Dehn fillings}
Let $M, M'$ be elements in $\calM_{g,k}$ with canonical decompositions
$\calT,\calT'$ respectively. Let $N$ (resp.~$N'$) be a hyperbolic
manifold obtained by Dehn filling $M$ (resp.~$M'$) along sufficiently
complicated slopes, and let $x\in I\Omega_{g,k}$ (resp.~$x'\in I\Omega_{g,k}$)
be such that $N\cong\widehat{M}(x)$, $N'\cong\widehat{M}' (x')$. 
In this paragraph we describe an explicit criterion which allows us to
determine if $N$ is commensurable with $N'$ just by looking at
$x,x'$ and at 
the combinatorics of $\calT,\calT'$.

\begin{defn}\label{commensurable:defn}
Two complete hyperbolic $n$-manifolds 
with geodesic boundary $M_1$, $M_2$ are
\emph{commensurable} if a hyperbolic manifold with geodesic boundary $M_3$ 
exists which is the total space of a finite Riemannian
covering both of $M_1$ and of $M_2$.
\end{defn}

\begin{prop}\label{comm:cond:prop}
Let $N_1,N_2$ be complete finite-volume hyperbolic $n$-manifolds
with non-empty geodesic boundary and denote by $\widetilde{N}_1,
\widetilde{N}_2$ the universal coverings of $N_1,N_2$ respectively.
Then $N_1$ is commensurable with $N_2$ if and only if
$\widetilde{N}_1$ is isometric to $\widetilde{N}_2$.
\end{prop}
\dimostraz
See~\cite{Fri:preprint2}.
\finedimo

From now on, let $k$ be a fixed \emph{odd} natural number
and let $X_k$ be the manifold defined in 
Subsection~\ref{nonhomeo:subsection}.

Let $\Delta_1,\ldots,\Delta_{2k+1}$ be the partially truncated tetrahedra 
of the canonical decomposition $\calT_k$ of $X_k$, and suppose as usual that 
for $i=1,\ldots,k$ the tetrahedra $\Delta_{2i-1},\Delta_{2i}$
are non-compact and glue up to a neighbourhood of the 
$i$-th cusp of $X_k$, while $\Delta_{2k+1}$ is compact regular.
We denote by
$F_{2i-1}^1,F_{2i-1}^2,F_{2i-1}^3,F_{2i}^1,
F_{2i}^2,F_{2i}^3$ the exceptional hexagons
of $\Delta_{2i-1},\Delta_{2i}$, in such a way that 
$F_{2i-1}^j$ is glued to $F_{2i}^j$
for $i=1,\ldots,k$, $j=1,2,3$. For $l=1,\ldots,2k$
we also call
$e_l^j$ the only finite edge of $F_l^j$, and $f_l^j$
the edge of $\Delta_l$ opposite to $e_l^j$.
We emphasize that here we do \emph{not} require that $F_l^1,F_l^2,F_l^3$
are positively arranged around the ideal vertex of $\Delta_l$.
Recall that a point $x\in\Omega_{k+1,k}$ determines
the geometric realization of $\calT_k$   
with dihedral angle $x_{6l-6+j}$ along $e_l^j$,
angle $x_{6l-3+j}$ along $f_l^j$, and angle
$x_{12k+1}$ along the compact edges of the unique compact tetrahedron.
It is easily seen that the exceptional lateral hexagons 
of the non-compact tetrahedra can be ordered around the ideal vertices
in such a way that the following condition
holds:
\begin{itemize}
\item
For $i=1,\ldots,k-1$, $j=1,2,3$ the isometry which glues
the compact face of $\Delta_{2i}$
to the compact face of $\Delta_{2i+1}$
sends $e_{2i}^j$ to $e_{2i+1}^j$.
Moreover, if $i$ is odd (resp.~even)
then $F^1_{2i-1},F^2_{2i-1},F^3_{2i-1}$ and 
$F^1_{2i},F^2_{2i},F^3_{2i}$ are positively (resp.~negatively)
arranged around the ideal
vertices of $\Delta_{2i-1}$ and $\Delta_{2i}$.  
\end{itemize}
(The fact that these conditions are coherent with each other depends
on the combinatorial properties of $\calT_k$. The second condition will
be taken into account when we will explicitly consider
the action of ${\rm Sym} (\Omega_{k+1,k})$ on $\Omega_{k+1,k}$).
Let $l\subset X_k$ be the compact edge of $\calT_k$. 
A straight-forward analysis of the combinatorics of $\calT_k$ 
shows that the dihedral angles of the tetrahedra of $\calT_k$
are arranged along $l$ according to the following cyclic ordering:
\[
\begin{array}{l}
x_{12k+1},\quad  x_1,x_7,x_{13}\ldots,x_{6l+1},\ldots, x_{12k-5},\\
x_{12k+1},x_{12k+1},\quad x_2,x_8,x_{14}\ldots,x_{6l+2},\ldots,x_{12k-4},\\
x_{12k+1},x_{12k+1},x_{12k+1},\quad x_3,x_9,x_{15}\ldots,x_{6l+3},
\ldots,x_{12k-3}.
\end{array}
\]
For $i=1,\ldots,k$ let $a_i,b_i,c_i:\Omega_{k+1,k}\to\matR$ 
be the functions defined as follows:
\[
\begin{array}{l}
a_i(x)=x_{12(j-1)+1}+x_{12(j-1)+7},\qquad 
b_i(x)=x_{12(j-1)+2}+x_{12(j-1)+8},\\
c_i(x)=x_{12(j-1)+3}+x_{12(j-1)+9}, 
\end{array}
\]
and set
\[
a,b,c:\Omega_{k+1,k}\to\matR,\qquad a(x)=\sum_{i=1}^{k} a_i(x),\quad
b(x)=\sum_{i=1}^{k} b_i(x),\quad
c(x)=\sum_{i=1}^{k} c_i (x).
\]
For $x\in\Omega_{k+1,k}$ we denote by $X_k (x)$
the hyperbolic structure defined on $X_k$ by $x$, and by
$\widehat{X}_k (x)$ the metric completion of 
$X_k (x)$. 
Let $V''$ be a
neighbourhood of $x_0$ in $\Omega_{k+1,k}$ as in the statement
of Proposition~\ref{univ:intrinsic:prop}  and for
$x\in V''\cap I\Omega_{k+1,k}$ let us denote by $U(x)$
the universal covering of $\widehat{X}_k (x)$.
We now show that the real numbers $a(x),b(x),c(x)$ completely determine
the isometry type of the universal covering
$U(x)$ of $\widehat{M} (x)$, whence the commensurability
class of $\widehat{M} (x)$.

\begin{prop}\label{comm:criterion:prop}
Let $x,x'$ be points in $I\Omega_{k+1,k}\cap V''$. Then 
$\widehat{X}_k (x)$ is commensurable with
$\widehat{X}_k (x')$ if and only if 
$a(x)=a(x'),\ b(x)=b(x'),\ c(x)=c(x')$.
\end{prop}
\dimostraz
Let $L^x (l_x)$ be the minimal distance between different connected components 
of $\partial U(x)$, and let  
$t\subset U(x)$ be a geodesic arc of length $L^x (l_x)$ 
joining two such components $S_1,S_2$.
By Corollary~\ref{short:return2:cor}, if $t'\subset U(x)$
is any other geodesic arc of length $L^x (l_x)$ connecting different
components of $\partial U(x)$, then there exists an isometry
of $U(x)$ taking $t$ to $t'$.
Let us consider the set $R\subset U(x)$ given by the union
of all the compact regular truncated tetrahedra whose truncation triangles
lie on $S_1\cup S_2\cup S'\cup S''$ for some
connected components $S',S''$ of $\partial U(x)$.
Let $N_\epsilon (t)$ be the $\epsilon$-neighbourhood of $t$, and consider
the sets $A=N_\epsilon (t)\cap R$ and $B=N_\epsilon (t)\setminus R$. 
Both $A$ and $B$ are unions of germs of dihedral sectors
whose number, amplitude and cyclic order (up to the choice of a positive
orientation around $t$) 
only depend on the isometry type of $U(x)$.
We will call such sectors $A$-sectors or $B$-sectors, according
to the fact that they
are contained
in $A$ or in $B$.
Lemma~\ref{short:return:lemma} implies that $t$ 
is a lift in $U(x)$ of the unique compact edge of the geometric
triangulation of $X_k (x)$, while by Proposition~\ref{univ:intrinsic:prop}
the set $R$ is the union of the lifts containing $t$ 
of the geometric tetrahedron $\Delta_{12k+1}\subset X_k (x)$.
Thus $A$-sectors are in number of three and have
angles $x_{12k+1},2x_{12k+1}$ and $3x_{12k+1}$. Moreover,
the $B$-sector between the $A$-sectors
with angles $x_{12k+1},2x_{12k+1}$ has angle
$a(x)$; the $B$-sector between the $A$-sectors
with angles $2x_{12k+1},3x_{12k+1}$ has angle
$b(x)$; the $B$-sector between the $A$-sectors
with angles $3x_{12k+1},x_{12k+1}$ has angle
$c(x)$.
This shows that $a(x),b(x),c(x)$ can be recovered solely from the
isometry type of $U(x)$, so if $\widehat{X}_k (x)$ is commensurable with
$\widehat{X}_k (x')$ we have $a(x)=a(x'), b(x)=b(x'),c(x)=c(x')$.

Suppose now that $a(x)=a(x'), b(x)=b(x'),c(x)=c(x')$.
Since $a(x)+b(x)+c(x)+6x_{12k+1}=a(x')+b(x')+c(x')+6x'_{12k+1}$
we have $x_{12k+1}=x'_{12k+1}$, so the compact tetrahedron in the
decomposition of $X_k (x)$ is isometric to the compact tetrahedron
in the decomposition of $X_k (x')$. Let now $U_\ast (x)$ (resp.~$U_\ast (x')$)
be the complement in $U(x)$ (resp.~in $U(x')$) of the preimage 
of the added geodesics
$\widehat{X_k} (x)\setminus X_k (x)$ (resp.~$\widehat{X_k} (x')\setminus X_k (x')$).
The geometric decomposition of $X_k$ parameterized by $x$ (resp.~$x'$)
naturally lifts to a tessellation $\calD_\ast (x)$ of $U_\ast (x)$
(resp.~$\calD_\ast (x')$ of $U_\ast (x')$). Let $K(x)\subset U(x)$ 
(resp.~$K(x')\subset U(x')$) be the union of the compact pieces of
$\calD_\ast (x)$ (resp.~$\calD_\ast (x')$).
Let now $U(x), U(x')$ be isometrically identified with
suitable polyhedra in $\matH^3$. 
Since the compact tetrahedra of $\calD_\ast (x)$ and $\calD_\ast (x')$
are all isometric to each other and $a(x)=a(x'), b(x)=b(x'),c(x)=c(x')$
it is easily seen that an element $\psi\in{\rm Isom}(\matH^3)$ 
exists which takes
$K(x)\subset U(x)\subset \matH^3$ onto $K(x')\subset U(x')\subset\matH^3$.
Since any component of $\partial U(x)$ (resp.~$\partial U(x')$) meets 
$K(x)$ (resp.~$\partial K(x')$) in a non-empty open subset of
a hyperbolic plane, this readily implies that $\psi (\partial U(x))=
\psi (\partial U(x'))$. Now $U(x),U(x')$ are the hyperbolic convex hulls
of $\partial U(x),\partial U(x')$ respectively, so $\psi (U(x))=U(x')$.
By Proposition~\ref{comm:cond:prop}, this implies that $\widehat{X}_k (x)$
and $\widehat{X}_k (x')$ are commensurable with each other.
\finedimo

In order to determine if geometrically similar manifolds
obtained by Dehn filling $X_k$ are commensurable with each other,
we are now reduced to understand when the functions $a,b$ and $c$
introduced above take different values on  
${\rm Sym} (\Omega_{k+1,k})$-equivalent points in
$\Omega_{k+1,k}$. 

Let us set
\[
H_h=\{x\in\matR^{12k+1}:\ x_{12i+1}=x_{12i+2}=x_{12i+3}
\ {\rm for\ all}\ i=h-1,\ldots,k-1\}.
\]
We recall that in a neighbourhood of $x_0$ in $\Omega_{k+1,k}$
the set $\Omega_{k+1,k}^h:=H_h\cap\Omega_{k+1,k}$ is a smooth manifold
of dimension $2h$ whose points
correspond to those structures which induce a complete metric
on the last $h$ cusps of $X_k$ (see Lemma~\ref{iso1:lemma}
and Proposition~\ref{iso4:prop}). 
Let now $\overline{\varsigma}_h:(-\varepsilon,\varepsilon)\to
\Omega^h_{k+1,k}$ be the curve mentioned in Remark~\ref{new:ext:rem}.
For a smooth $f:\Omega_{k+1,k}\to\matR$ let us denote by
$\dot{f}$ (resp.~by $\ddot{f}$) the first (resp.~second) 
derivative of $f\circ\overline{\varsigma}_h$ at $0$.
From Proposition~\ref{fine:der:prop} we deduce:
\[
\begin{array}{llllll}
\dot{a}_i&=&\dot b_i&=&\dot c_i=0 & {\rm for\ all}\ i=1,\ldots,k;\\
\ddot{a}_i&>&\ddot b_i&>&\ddot c_i & {\rm for\ all}\ i=1,\ldots,h;\\
\ddot{a}_i&=&\ddot b_i&=&\ddot c_i=0 & {\rm for\ all}\ i=h+1,\ldots,k.
\end{array}
\]

We are now ready to prove the following:

\begin{teo}\label{noncomm:fil:teo}
Fix $1\leqslant h\leqslant k$, where $k$ is odd. Then there exists a sequence
$\{W_h^n\}_{n\in\matN}$ of pairwise non-homeomorphic complete finite-volume
hyperbolic manifolds with geodesic boundary with the following properties:
\begin{itemize}
\item
Each $W_h^n$ is obtained by Dehn filling the first $h$
cusps of $X_k$;
\item
For any $n\in\matN$
there exist at least three (including $W_h^n$ itself)
pairwise non-commensurable
hyperbolic Dehn fillings of $X_k$ 
which are geometrically similar to $W_h^n$.
\end{itemize}
\end{teo}
\dimostraz
We choose
an infinite sequence $\{y_n\}_{n\in\matN}\subset 
I\Omega_{k+1,k}^h\setminus \{x_0\}$
converging to $x_0$
along  
$\dot{\overline{\varsigma}}_h (0)$ (see Definition~\ref{direction:defn}), 
and we set $W_h^n=\widehat{X}_k (y_n)$.

We first observe that $W_h^n$ is obtained from $X_k$ by Dehn filling
the first $h$ cusps of $X_k$: since $y_n$ belongs
to $\Omega_{k+1,k}^h$, the last $k-h$ cusps of $X_k (y_n)$ have to be complete;
moreover, up to extracting a subsequence
we can suppose $a_i (y_n)>b_i (y_n)>c_i (y_n)$ for all $i=1,\ldots,h$,
so the angles along the compact edges of $\Delta_{2i-1},\Delta_{2i}$
are not equal to each other, and the $i$-th cusp of
$X_k (y_n)$ is not complete.

Let now $r\in\calM (T_1\sqcup\ldots\sqcup T_k)$
be the element acting as a 
positive (resp.~negative) rotation
by an angle of
$\pi/3$ on $T_i$ for $i\in\{1,\ldots,k\}$, $i$ odd (resp.~even),
and let $\Theta:\calM(T_1\sqcup\ldots\sqcup T_k)\to
{\rm Sym} (\Omega_{k+1,k})$ be the isomorphism described in
Proposition~\ref{iso:to:iso:prop}. We set
$y'_n=\Theta (r) (y_n)$ and $y''_n=\Theta (r^2) (y_n)$.
By construction,  $W_h^n=\widehat{X}_k (y_n),
\widehat{X}_k (y'_n)$ and
$\widehat{X}_k (y''_n)$ are pairwise geometrically similar.
Moreover an easy computation shows that for
$x\in\Omega_{k+1,k}$ we have
\[
\begin{array}{c}
a_i (\Theta (r^2) (x))=c_i (\Theta (r)(x))=b_i (x),\\
b_i (\Theta (r^2) (x))=a_i (\Theta (r)(x))=c_i (x),\\
c_i (\Theta (r^2) (x))=b_i (\Theta (r)(x))=a_i (x),
\end{array}
\]
whence $a_i (y_n)>a_i (y''_n)>a_i (y'_n)$,
and $W_h^n=\widehat{X}_k (y_n),
\widehat{X}_k (y'_n), 
\widehat{X}_k (y''_n)$ are pairwise
non-commensurable by 
Proposition~\ref{comm:criterion:prop}.
\finedimo

\begin{rem}
Let $M$ be an element of $\calM_{g,k}$ with canonical decomposition
$\calT$. Suppose that 
the arrangement of compact and non-compact
tetrahedra around the
compact edge of $\calT$ is sufficiently irregular
and let $\calS=\{s_{i_1},
\ldots,s_{i_h}\}$ be 
a set of slopes for $M$ such that $s_{i_l}$ is not equivalent
to $s_{i_{m}}$ for $l\neq m$.  
The same argument used to prove 
Theorem~\ref{noncomm:fil:teo}
shows that the Dehn fillings of $M$
which are geometrically similar to $M(\calS)$ are expected to be
non-commensurable with each other.
\end{rem}

We conclude with some examples of non-homeomorphic
geometrically similar commensurable Dehn fillings of $X_k$.

\begin{teo}\label{comm:fll:teo}
Let $k\geqslant 3$ be odd. Then there exists an infinite sequence
of pairs $\{Y_1^n, Y_2^n \}_{n\in\matN}$ of complete finite-volume
hyperbolic manifolds with geodesic boundary such that for every
$n\in\matN$ the following conditions hold:
\begin{itemize}
\item
$Y_1^n$ is obtained by Dehn filling the first
cusp of $X_k$;
\item
$Y_2^n$ is obtained by Dehn filling the third
cusp of $X_k$;
\item
$Y_1^n$ is geometrically similar to $Y_2^n$;
\item
$Y_1^n$ is commensurable with $Y_2^n$;
\item
$Y_1^n$ is not homeomorphic to $Y_2^n$. 
\end{itemize}
\end{teo}
\dimostraz
We choose
an infinite sequence $\{y_n\}_{n\in\matN}\subset
I\Omega_{k+1,k}^1\setminus\{x_0\}$
converging to $x_0$ along
$\dot{\overline{\varsigma}} (0)$.
Let $\tau_{13}\in{\rm Sym} (\Omega_{k+1,k})$ 
be the element which exchanges the first
cusp of $X_k$ with the third one according to equation~(\ref{def:kappa:eq})
and let 
$y'_n=\tau_{13} (y_n)$. We set $Y_1^n=\widehat{X}_k (y_n)$
and $Y_2^n=\widehat{X}_k (y'_n)$.
It is easily seen that $Y_1^n$ is obtained by filling the
first cusp of $X_k$, while $Y_2^n$ is obtained by filling 
the third one.
The element of $\calM (T_1\sqcup\ldots\sqcup T_k)$ corresponding
to $\tau_{13}$ is orientation-preserving, so
$Y_1^n$ is geometrically similar to $Y_2^n$.
Moreover an easy computation shows that for every
$x\in\Omega_{k+1,k}$ we have
\[
\begin{array}{llll}
a_1 (\tau_{13} (x))=a_3 (x),&
b_1 (\tau_{13} (x))=b_3 (x),&
c_1 (\tau_{13} (x))=c_3 (x),& \\
a_3 (\tau_{13} (x))=a_1 (x),&
b_3 (\tau_{13} (x))=b_1 (x),&
c_3 (\tau_{13} (x))=c_1 (x),& \\
a_j (\tau_{13} (x))=a_j (x),&
b_j (\tau_{13} (x))=b_j (x),&
c_j (\tau_{13} (x))=c_j (x),&
j=2,4,5,6,\ldots,k.
\end{array}
\]
This easily implies
$a (y_n)=a (y'_n)$, 
$b (y_n)=b (y'_n)$, $c (y_n)=c (y'_n)$,
so $Y_1^n$ is 
commensurable with $Y_2^n$ by 
Proposition~\ref{comm:criterion:prop}.

Let us prove that $Y_1^n$ is not homeomorphic to $Y_2^n$.
Up to passing to a subsequence we can suppose that the added
geodesic $Y_1^n\setminus X_k$ (resp.  $Y_2^n\setminus X_k$)
is the shortest geodesic of $Y_1^n$ (resp. $Y_2^n$). 
Let $f_n:Y_1^n \to Y_2^n$ be a homeomorphism. By Mostow-Prasad's
rigidity Theorem we may assume  that $f_n$ is an isometry, 
which implies $f_n (Y_1^n\setminus X_k)=Y_2^n\setminus X_k$.
Thus $f_n$ restricts to a homeomorphism $f'_n: X_k\to X_k$. 
By rigidity again 
we can homotope $f'_n$ into an isometry, which by construction should
take the first cusp of $X_k$ onto the third one, against 
Proposition~\ref{iso:xk:prop}.
\finedimo

\section*{Acknowledgements}
The results described in this paper are contained in my doctoral thesis.
I would like to thenk Prof.~Carlo Petronio, my thesis advisor, for his
support and encouragement.

\bibliographystyle{amsalpha}
\bibliography{bibliosim}

\end{document}